\documentclass[12pt]{amsart}
\usepackage{amsmath,amssymb,amsbsy,amsfonts,amsthm,latexsym, mathtools,amsopn,amstext,amsxtra,euscript,amscd,mathrsfs,color,bm, cite,multirow,tabularx}
\usepackage{graphicx} 
\usepackage{todonotes,dirtytalk}   
\usepackage{url}
\usepackage[dvipsnames]{xcolor}
\usepackage[colorlinks,linkcolor=blue,anchorcolor=blue,citecolor=blue]{hyperref}
\usepackage{color}
\usepackage{comment}
\allowdisplaybreaks

\DeclareMathOperator{\li}{li}

\begin{document}
\newtheorem{problem}{Problem}
\newtheorem{theorem}{Theorem}
\newtheorem{lemma}[theorem]{Lemma}
\newtheorem{crit}[theorem]{Criterion}
\newtheorem{claim}[theorem]{Claim}
\newtheorem{cor}[theorem]{Corollary}
\newtheorem{prop}[theorem]{Proposition}
\newtheorem{definition}{Definition}
\newtheorem{question}[theorem]{Question}
\newtheorem{Note}[theorem]{Notation}
\newtheorem{hyp}{Hypothesis}

\theoremstyle{remark} 
\newtheorem*{note}{Note}
\newtheorem*{rem}{Remark}

\def\cA{{\mathcal A}}
\def\cB{{\mathcal B}}
\def\cC{{\mathcal C}}
\def\cD{{\mathcal D}}
\def\cE{{\mathcal E}}
\def\cF{{\mathcal F}}
\def\cG{{\mathcal G}}
\def\cH{{\mathcal H}}
\def\cI{{\mathcal I}}
\def\cJ{{\mathcal J}}
\def\cK{{\mathcal K}}
\def\cL{{\mathcal L}}
\def\cM{{\mathcal M}}
\def\cN{{\mathcal N}}
\def\cO{{\mathcal O}}
\def\cP{{\mathcal P}}
\def\cQ{{\mathcal Q}}
\def\cR{{\mathcal R}}
\def\cS{{\mathcal S}}
\def\cT{{\mathcal T}}
\def\cU{{\mathcal U}}
\def\cV{{\mathcal V}}
\def\cW{{\mathcal W}}
\def\cX{{\mathcal X}}
\def\cY{{\mathcal Y}}
\def\cZ{{\mathcal Z}}

\def\A{{\mathbb A}}
\def\B{{\mathbb B}}
\def\C{{\mathbb C}}
\def\D{{\mathbb D}}
\def\E{{\mathbb E}}
\def\F{{\mathbb F}}
\def\G{{\mathbb G}}
\def\I{{\mathbb I}}
\def\J{{\mathbb J}}
\def\K{{\mathbb K}}
\def\L{{\mathbb L}}
\def\M{{\mathbb M}}
\def\N{{\mathbb N}}
\def\O{{\mathbb O}}
\def\P{{\mathbb P}}
\def\Q{{\mathbb Q}}
\def\R{{\mathbb R}}
\def\S{{\mathbb S}}
\def\T{{\mathbb T}}
\def\U{{\mathbb U}}
\def\V{{\mathbb V}}
\def\W{{\mathbb W}}
\def\X{{\mathbb X}}
\def\Y{{\mathbb Y}}
\def\Z{{\mathbb Z}}

\def\ep{{\mathbf{e}}_p}
\def\eq{{\mathbf{e}}_q}
\def\cal#1{\mathcal{#1}}

\def\scr{\scriptstyle}
\def\\{\cr}
\def\({\left(}
\def\){\right)}
\def\[{\left[}
\def\]{\right]}
\def\<{\langle}
\def\>{\rangle}
\def\fl#1{\left\lfloor#1\right\rfloor}
\def\rf#1{\left\lceil#1\right\rceil}
\def\le{\leqslant}
\def\ge{\geqslant}
\def\eps{\varepsilon}
\def\mand{\qquad\mbox{and}\qquad}

\def\sssum{\mathop{\sum\ \sum\ \sum}}
\def\ssum{\mathop{\sum\, \sum}}
\def\ssumw{\mathop{\sum\qquad \sum}}

\def\vec#1{\mathbf{#1}}
\def\inv#1{\overline{#1}}
\def\num#1{\mathrm{num}(#1)}
\def\dist{\mathrm{dist}}

\def\fA{{\mathfrak A}}
\def\fB{{\mathfrak B}}
\def\fC{{\mathfrak C}}
\def\fU{{\mathfrak U}}
\def\fV{{\mathfrak V}}

\newcommand{\bflambda}{{\boldsymbol{\lambda}}}
\newcommand{\bfxi}{{\boldsymbol{\xi}}}
\newcommand{\bfrho}{{\boldsymbol{\rho}}}
\newcommand{\bfnu}{{\boldsymbol{\nu}}}

\def\GL{\mathrm{GL}}
\def\SL{\mathrm{SL}}
\def\Ei{\mathrm{Ei}}

\def\Hba{\overline{\cH}_{a,m}}
\def\Hta{\widetilde{\cH}_{a,m}}
\def\Hb1{\overline{\cH}_{m}}
\def\Ht1{\widetilde{\cH}_{m}}

\def\flp#1{{\left\langle#1\right\rangle}_p}
\def\flm#1{{\left\langle#1\right\rangle}_m}

\def\Zm{\Z/m\Z}

\def\Err{{\mathbf{E}}}
\def\O{\mathcal{O}}

\def\cc#1{\textcolor{red}{#1}}
\newcommand{\commT}[2][]{\todo[#1,color=green!60]{Tim: #2}}
\newcommand{\commB}[2][]{\todo[#1,color=red!60]{Bryce: #2}}

\DeclarePairedDelimiter\ceil{\lceil}{\rceil}
\DeclarePairedDelimiter\floor{\lfloor}{\rfloor}

\newcommand{\comm}[1]{\marginpar{%
\vskip-\baselineskip 
\raggedright\footnotesize
\itshape\hrule\smallskip#1\par\smallskip\hrule}}
\newcolumntype{L}{>{\raggedright\arraybackslash}X}

\def\xxx{\vskip5pt\hrule\vskip5pt}

\def\dmod#1{\,\left(\textnormal{mod }{#1}\right)}


\title[Explicit Burgess bounds for prime moduli]{Improvements to explicit Burgess-strength bounds for quadratic characters of prime modulus}

\author[J. Punch]{James Punch}\thanks{}
\address{School of Science, The University of New South Wales Canberra, Australia}
\email{j.punch@unsw.edu.au}

\pagenumbering{arabic}

\begin{abstract}
For positive integers $M$ and $N$, consider the character sum $S(M,N)=\sum_{n\in[M,M+N)}\chi_p(n)$, where $\chi_p$ is the quadratic character modulo $p$. Burgess \cite{Burgess1963} proved that $S(M,N)\ll N^{1-1/r}p^{\frac{r+1}{4r^2}+\varepsilon}$ for any integer $r\geq 2$ and any real $\varepsilon>0$. This was made explicit by Booker \cite{Booker2006}, who showed (under size restrictions on $M$, $N$, $p$ and $r$) that  $\lvert S(M,N)\rvert \leq \alpha(r)N^{1-1/r}p^{\frac{r+1}{4r^2}}(\beta(r)+\log{p})^{1/2r}$, where $\alpha(r)$ and $\beta(r)$ are given in \cite{Booker2006}. We show that when $M$ is not too small compared to $N$, we can improve $\beta(r)$ by a fixed amount. If $M$ is also not too large compared to $N$, we can further improve $\beta(r)$ by an amount proportional to $\frac{M}{M+N}$.
\end{abstract}

\maketitle
\let\thefootnote\relax
\footnote{\textit{Affiliation}: School of Science, University of New South Wales, Canberra.}

\section{Introduction}
Let $p$ be an odd prime. For a non-principal character $\chi$ with modulus $p$, explicit upper bounds for the character sum $S(M,N)=\sum_{n\in[M,M+N)}\chi(n)$ are regularly required in number-theoretic settings, such as bounding the least quadratic non-residue modulo $p$ (see the survey article by \cite{McGownTrevino2021}) or the least primitive root modulo $p$ (see \cite{Grosswald1981} for a conjecture of Grosswald and \cite{CohenOliveiraeSilvaTrudgian2016}, \cite{McGownTrevinoTrudgian2016}, \cite{McGownTrudgian2020}, \cite{McGownSorenson2025} and \cite{Bagger2026} for recent progress). They can also be used to compute values of Dirichlet $L$-functions, which may, via the class number formula, be used to compute class numbers of real quadratic fields (see \cite{Booker2006}) and ranks of elliptic curves (see \cite{Elkies2025}).

It was proved by Burgess \cite{Burgess1963} that $S(M,N)\ll N^{1-1/r}p^\frac{r+1}{4r^2}\log{p}$, and the power of $\log{p}$ was improved from $1$ to (in particular) $\frac{1}{2r}$ by Montgomery and Vaughan \cite{MV2007}. We note that Iwaniec and Kowalski \cite{IwaniecKowalski2004} gave $\frac{1}{r}$ before this and that other authors have since improved it; the most recent result is $0$, proved by Kerr, Trudgian and the author \cite{KerrPunchTrudgian2026}.

From now on, unless specified otherwise, $S(M,N)$ will always mean a character sum involving the quadratic character modulo $p$. The result of Montgomery--Vaughan was made explicit for quadratic characters by Booker \cite{Booker2006}, who proved the following.

\begin{theorem} \label{prop:Booker}
    Let $p>10^{20}$ be a prime number equivalent to $1\pmod{4}$, select $r\in\{2,3,\ldots,15\}$ and choose $M,N$ to be positive integers not exceeding $2\sqrt{p}$. Then $$\lvert S(M,N)\rvert \leq \alpha_B(r) N^{1-1/r}p^{\frac{r+1}{4r^2}}(\beta_B(r)+\log{p})^{1/2r},$$ where $\alpha_B(r)$ and $\beta_B(r)$ are given in Table \ref{tab:alphabeta}.
\end{theorem}

In most Burgess-strength results (including \cite{Booker2006}), one considers the function $\nu(x)$, which counts the number of ways to write some equivalence class $x\pmod{p}$ as $a^{-1}n$, where $\lvert a\rvert$ is an integer less than or equal to some parameter $A$, and $n\in [M,M+N)$. Specifically, an application of the H\"older inequality leads us to bound the first and second moments of $\nu(x)$ (i.e.\ $\sum_{x=1}^{p}\nu(x)$ and $\sum_{x=1}^{p}\nu(x)^2$). This paper originated from an investigation into whether leveraging higher moments can lead to an improvement; we will prove that it cannot. However, this work also led to an improved bound on the second moment of $\nu(x)$. This bound is unusual because it improves as $M$ becomes large in comparison to $N$; the author is not aware of any other result of this type.

Now, Burgess-style proofs involve shifting and averaging the desired character sum. \cite{Booker2006} then uses proof by induction (in the same style as \cite{IwaniecKowalski2004}), but \cite{MV2007} adopts a recursive approach instead. We show that it is possible to use the dependence on $M$ to our advantage in both styles of the proof. We state the Booker-style result here and prove it in Section \ref{sec:BurgessBooker}, but defer the statement of the Montgomery--Vaughan-style result (which has worse constants but is more flexible) and its proof to Section \ref{sec:BurgessMV}.

\begin{table}
    \centering
    \begin{tabular}{|ccccc|ccccc|}
    \hline
        $r$ & $\alpha(r)$ & $\alpha_B(r)$ & $\beta(r)$ & $\beta_B(r)$ & $r$ & $\alpha(r)$ & $\alpha_B(r)$ & $\beta(r)$ & $\beta_B(r)$\\ \hline
        2 & 1.8221 & 1.8221 & -7.1055 & 8.9077 & 9 & 1.4091 & 1.4548 & -8.9989 & 0.0085\\
        3 & 1.7157 & 1.8000 & -6.6151 & 5.3948 & 10 & 1.3819 & 1.4231 & -9.3068 & -0.4106\\
        4 & 1.6410 & 1.7263 & -7.0097 & 3.6658 & 11 & 1.3584 & 1.3958 & -9.5921 & -0.7848\\
        5 & 1.5767 & 1.6526 & -7.4677 & 2.5405 & 12 & 1.3378 & 1.3721 & -9.8577 & -1.1232\\
        6 & 1.5229 & 1.5892 & -7.9020 & 1.7059 & 13 & 1.3197 & 1.3512 & -10.1061 & -1.4323\\
        7 & 1.4782 & 1.5363 & -8.3005 & 1.0405 & 14 & 1.3037 & 1.3328 & -10.3394 & -1.7169\\
        8 & 1.4408 & 1.4921 & -8.6648 & 0.4856 & 15 & 1.2893 & 1.3164 & -10.5593 & -1.9808\\ \hline
    \end{tabular}
    \caption{Table of $\alpha(r)$, $\alpha_B(r)$, $\beta(r)$ and $\beta_B(r)$ as in Theorem \ref{prop:BurgessBooker}}
    \label{tab:alphabeta}
\end{table}
\begin{theorem} \label{prop:BurgessBooker}
    Let $p$ be a prime number greater than $10^{23}$; let $M_0$ and $N_0$ be positive integers with $N_0\leq 2\sqrt{p}$ and select some integer $r\in\{2,3,\ldots, 15\}$. Write $c_0=\frac{M_0}{M_0+N_0}$ and let $c_0'=\max(1-c_0,13/61)$. Let $k=\left(\frac{2r-1}{2r^2}\right)^{r/(r-1)}$ and suppose that $c_0\geq k$. Further, let $C_0=\frac{p^{\frac{1}{2}+\frac{1}{2r}}}{8k}\left[\frac{1}{2}(2r-1)!!\right]^{\frac{1}{r}}$ and suppose that $M_0+\frac{N_0}{1-k}\leq C_0$. Let $$K=\alpha(r)\left(\log{p}+\beta(r)+\frac{4r}{r-1}c_0'\right)^{1/2r},$$ where $\alpha(r)$ and $\beta(r)$ are given in Table \ref{tab:alphabeta}. Then $\lvert S(M_0,N_0)\rvert \leq KN_0^{1-1/r}p^{\frac{r+1}{4r^2}}$. 
\end{theorem}
\begin{rem}
    Except when $r=2$, each $\alpha(r)$ is smaller than $\alpha_B(r)$. The author does not claim this as an improvement, because it actually arises from correcting what he believes is a mistake in \cite{Booker2006}. We will give more detail about this later.
\end{rem}

Finally, we note that computation is required at many points in this paper. All computations were done in \emph{Mathematica} and they can be seen at \url{https://github.com/JP1681/Burgess-bounds}.

\section{Behaviour of higher moments of \texorpdfstring{$\nu(x)$}{ν(x)}} \label{sec:higherMoments}
Let $\mathcal{A}$ and $\mathcal{N}$ be sets of integers, and let $\nu(x)=\nu_{\mathcal{A},\mathcal{N}}(x)$ be the number of representations of $x$ as $a^{-1}n\pmod{p}$ where $a\in \mathcal{A}$ and $n\in \mathcal{N}$. A crucial step in the argument of \cite{Booker2006} involves the following expression:
\begin{align}\label{eq:Booker23}\sum_{x\pmod{p}}\nu_{\mathcal{A},\mathcal{N}}(x)\left\lvert\sum_{b=-B}^{B} \chi_p(x+b)\right\rvert,\end{align} where $\mathcal{A}$ is the set of primes whose magnitudes do not exceed some chosen value $A<p$, and $\mathcal{N}=[M,M+N)\,\cap\,\Z$ for integers $M$ and $N$.

Booker uses H\"older's inequality with the identity $\nu(x)=[\nu(x)]^{1-1/r}[\nu(x)^2]^{1/2r}$ to reach an upper bound of: $$\left[\sum_{x\pmod{p}}\nu(x)\right]^{1-1/r}\left[\sum_{x\pmod{p}}\nu(x)^2\right]^{1/2r}\left[\sum_{x\pmod{p}}\left\lvert\sum_{b}\chi_p(x+b)\right\rvert^{2r}\right]^{1/2r}.$$ The first of the three sums over $x\pmod{p}$ is $\#\mathcal{A}\times \# \mathcal{N}=2\pi(A)N$. Bounding the second and third sums requires more effort.

This strategy naturally leads to the question: why stop at the second moment of $\nu(x)$? We now generalise: if $m\geq 2$ and $a_1+2a_2+\ldots+ ma_m=1$, then the above identity becomes $\nu(x)=[\nu(x)]^{a_1}[\nu(x)^2]^{a_2}\ldots [\nu(x)^m]^{a_m}$. Assuming further that $a_1+\ldots + a_m+\frac{1}{2r}=1$, H\"older's inequality gives an upper bound of \begin{align}\label{eq:HolderBound}\left[\sum_{x\pmod{p}}\nu(x)\right]^{a_1}\dots\left[\sum_{x\pmod{p}}\nu(x)^m\right]^{a_m}\left[\sum_{x\pmod{p}}\left\lvert\sum_{b}\chi_p(x+b)\right\rvert^{2r}\right]^{1/2r}.\end{align} Taking $m=2$ forces $a_1=1-1/r$ and $a_2=1/2r$ as above.

Booker shifts by integers $ab$, which may be positive, negative, or zero; we will do the same in Section \ref{sec:BurgessBooker}. Our proof in Section \ref{sec:BurgessMV} will employ shifts by positive integers only, so we must consider \begin{align}\label{eq:Booker23modified}\sum_{x\pmod{p}}\nu_{\mathcal{A},\mathcal{N}}(x)\left\lvert\sum_{b=1}^{B} \chi_p(x+b)\right\rvert,\end{align} where $\mathcal{A}$ is now \emph{the set of primes not exceeding $A$}, and $\mathcal{N}$ is as before.
Using H\"older as above, we would need to bound $\sum_{x\pmod{p}}\nu(x)^k$ for each $2\leq k\leq m$. We now give a method for doing so that establishes a dependence on $M$, which we will use to our advantage later.

\subsection{An expression for the \texorpdfstring{$k^\text{th}$}{k-th} moment of \texorpdfstring{$\nu(x)$}{ν(x)}}
Let $\nu(x)=\nu_{\mathcal{A},\mathcal{N}}(x)$, where $\mathcal{A}$ is the set of primes not exceeding $A$ and $\mathcal{N}=[M,M+N)\cap \Z$. Our objective here is to derive an expression for $\sum_{x\pmod{p}}\nu(x)^k$ in terms of $M,N$ and $A$. In \cite{Booker2006}, the expression \eqref{eq:Booker23} is bounded under the assumption that $A(M+N)\leq p/2$. We deal with \eqref{eq:Booker23modified} by instead assuming that $A(M+N)\leq p$.

In what follows, all equivalences are modulo $p$, and the $a_i$ are always prime.

\begin{lemma} \label{prop:exactKthMoment}
    Let $p$ be a prime; suppose that $A, M, N$ and $k$ are positive integers and that $A(M+N)\leq p$. Write $P=P(a_1,a_2):=\pi(a_1)-\pi(a_2)$. If $k\geq 2$, then \begin{align*}
    \sum_{x\pmod{p}}\nu(x)^k&=\pi(A)N +\\&\sum_{a_1\leq A}\sum_{\frac{M}{M+N}a_1<a_2<a_1}\left(\ceil*{\frac{M+N}{a_1}}-\ceil*{\frac{M}{a_2}}\right)\left((P+1)^k-2P^k+(P-1)^k\right).
\end{align*}
\end{lemma}
\begin{proof}
We have \begin{align*}
    \sum_{x\pmod{p}}\nu(x)^k&=\sum_{x\pmod{p}}\left[\sum_{a\in \mathcal{A}} \sum_{\substack{n\in \mathcal{N}\\a^{-1}n\equiv x\,}}1\right]^k=\sum_{x\pmod{p}}\left[\sum_{a_1,\ldots,a_k\in \mathcal{A}}\sum_{\substack{n_1,\ldots,n_k\in\mathcal{N}\\\forall i,a_i^{-1}n_i\equiv x}}1\right]\\
    &=\sum_{a_1,\ldots,a_k\in\mathcal{A}}\sum_{n_1,\ldots,n_k\in\mathcal{N}}\sum_{\substack{x\pmod{p}\\\forall i,x\equiv a_i^{-1}n_i}}1.
\end{align*}
Since $A(M+N)\leq p$, each $a_i n_j$ is in the interval $(0,p)$. Thus
\begin{align*}
    x\equiv a_i^{-1}n_i \text{ for all } i &\iff a_i n_j\equiv a_j n_i  \text{ for all } i,j\text{ and } a_1^{-1}n_1\equiv x\,\\
    &\iff a_i n_j= a_j n_i \text{ for all } i,j\text{ and } a_1^{-1}n_1\equiv x\, \\
    &\iff \frac{n_i}{a_i}= \frac{n_j}{a_j} \text{ for all } i,j \text{ and } a_1^{-1}n_1\equiv x\, .
\end{align*}
In particular, $\frac{n_i}{a_i}$ takes the same value for every $i$. The above sum is then \begin{align}\label{eq:niai}\sum_{a_1,\ldots,a_k\in \mathcal{A}}\sum_{\substack{n_1,\ldots,n_k \in \mathcal{N}\\\frac{n_i}{a_i}\text{ constant}}}\sum_{\substack{\\x\equiv a_1^{-1}n_1}}1=\sum_{a_1,\ldots,a_k \in \mathcal{A}}\sum_{\substack{n_1,\ldots,n_k \in \mathcal{N}\\\frac{n_i}{a_i}\text{ constant}}}1.\end{align}
If the $a_i$ are all the same, so too are the $n_i$. Since there are $\pi(A)$ options for $a_1$ and $N$ options for $n_1$, the contribution to the sum (the \say{diagonal term}) is $\pi(A)N$.

Otherwise, the $a_i$ are not all the same. We may then find $r,s$ for which $a_r \neq a_s$, and we must have $n_s a_r=n_r a_s$. Then $a_r\mid n_r$ and $a_s\mid n_s$ by Euclid's lemma; so $\frac{n_r}{a_r}=\frac{n_s}{a_s}\in \Z$. We now count the number of ways to make each of the $\frac{n_i}{a_i}$ equal to the same integer.

It turns out that this number does not depend on all the $a_i$ at once; only on $\min_i a_i$ and $\max_i a_i$. We therefore fix these two values and find the contribution to the sum, then count the ways that the minimum and maximum can take those two values.

Let $x=\min_i a_i$ and $y=\max_{i} a_i$, so that $x<y$. We have $\frac{n_i}{a_i}\in \left[\frac{M}{a_i},\frac{M+N}{a_i}\right)$ for each $i$, so we need to count the number of integers in the set $\left[\frac{M}{a_1},\frac{M+N}{a_1}\right) \cap \ldots \cap \left[\frac{M}{a_k},\frac{M+N}{a_k}\right).$ If $\frac{M+N}{y}\leq \frac{M}{x}$, this is the empty set. Put another way, we get no contribution if $x\leq \frac{M}{M+N}y$. Otherwise, $\frac{M}{M+N}y<x<y$, and the above set is the half-open interval $\left[\frac{M}{x},\frac{M+N}{y}\right)$, which contains $\ceil{(M+N)/y}-\ceil{M/x}$ integers.

We now need to know how many tuples $(a_1,\ldots,a_k)$ have their smallest and largest elements equal to $x$ and $y$ respectively. This is not difficult; first, let $X=\pi(x)$ and $Y=\pi(y)$, then each of the $a_i$ must be one of the $Y-X+1$ primes between $x$ and $y$ inclusive: $(Y-X+1)^k$ possibilities. However, we must use the $Y^\text{th}$ prime at least once, so subtract $(Y-X)^k$. We must also use the $X^\text{th}$ prime at least once, so remove another $(Y-X)^k$. The number of cases using neither the $X^\text{th}$ nor $Y^\text{th}$ primes was just subtracted twice, so add $(Y-X-1)^k$. Using $P$ to denote $\pi(y)-\pi(x)$, the number of ways is $(P+1)^k-2P^k+(P-1)^k$. Combining our results yields the lemma.
\end{proof}

\begin{rem}
    If we revert to the $\mathcal{A}$ from \cite{Booker2006} (where $-a$ and $a$ are included) and assume that $A(M+N)\leq p/2$, then the bound from Lemma \ref{prop:exactKthMoment} is doubled. The new assumption forces each $a_i n_j$ to be in the interval $(-p/2,p/2)$, so we may still say that $a_in_j\equiv a_j n_i \iff a_in_j=a_j n_i$. We double the bound since each solution in \eqref{eq:niai} where all $a_i$ are positive gives rise to exactly one other (obtained by multiplying each $a_i$ by $-1$).
\end{rem}

\subsection{An upper bound when \texorpdfstring{$k=2$}{k=2}}
We will later show that higher moments of $\nu(x)$ (i.e.\ $k>2$) cannot actually give any improvement. Fortunately, using our Lemma \ref{prop:exactKthMoment} on the second moment will lead to an improvement, so we now give an explicit upper bound for it. Again, we assume that $a_1$ and $a_2$ are always prime.

\begin{lemma} \label{prop:2ndMomentUpperBound}
    Let $p$ be a prime; suppose that $A$, $M$ and $N$ are positive integers, that $A\geq 3$ and that $A(M+N)\leq p$. Write $c=\frac{M}{M+N}$ and $c'=\max(1-c,\frac{13}{61})$. Then 
    $$\sum_{x\pmod{p}}\nu(x)^2\leq \pi(A)N+2N\left[\frac{122}{105}+c'\sum_{5\leq i \leq \pi(A)}\frac{i}{p_i}\right]+\left[12+c'(\pi(A)^2+\pi(A)-16)\right].$$
\end{lemma}
\begin{rem}
    In the following proof, we split the off-diagonal term from Lemma \ref{prop:exactKthMoment} into two components, and bound each separately. We refer to them as a \say{main term} and an \say{error term}, but this characterisation is only accurate if $A$ is much smaller than $N$. This lemma will later be applied in circumstances where $A/N$ is indeed very small, but if we fixed $N$ and made $A$ large, then the \say{error term} would dominate the \say{main term}.
\end{rem}
\begin{proof} We employ Lemma \ref{prop:exactKthMoment}. When $k=2$, we have $(P+1)^k-2P^k+(P-1)^k=2$. Forgetting for now the diagonal contribution $\pi(A)N$, the remaining part has main term
$$2\sum_{a_1\leq A}\sum_{ca_1<a_2<a_1}\left(\frac{M+N}{a_1}-\frac{M}{a_2}\right),$$ plus an error term whose magnitude is no larger than $2\sum_{a_1\leq A}\sum_{ca_1<a_2<a_1}1.$
Observe that $$\frac{M+N}{a_1}-\frac{M}{a_2}=\frac{N}{a_1}-M\left(\frac{1}{a_2}-\frac{1}{a_1}\right)< \frac{N}{a_1}.$$ The main term is then no greater than
\begin{align*}
2N\sum_{a_1\leq A}\frac{1}{a_1}\sum_{ca_1<a_2<a_1}1&= 2N\sum_{2\leq i\leq \pi(A)}\frac{1}{p_i}(i-\pi(cp_i)-1)\\
&\leq 2N\left[\frac{1}{3}+\frac{2}{5}+\frac{3}{7}+c'\sum_{5\leq i\leq \pi(A)}\frac{i}{p_i}\right],\end{align*} via Corollary \ref{prop:PunchConjectureUniform} with $n=1$.

If $\pi(A)\geq 5$, the size of the error term does not exceed \begin{align*}
    2\sum_{2\leq i\leq \pi(A)}(i-\pi(cp_i)-1)&\leq 2\left[6+c'\sum_{5\leq i \leq \pi(A)}i\right]\\
    &= 2\left[6+c'\left(\frac{\pi(A)(\pi(A)+1)}{2}-10\right)\right].
\end{align*} By computation, we can see that the bound also holds for $2\leq \pi(A)< 5$ if we replace the $-10$ by $-8$. Combining our results yields the lemma.
\end{proof}

\subsection{What happens when \texorpdfstring{$k\geq 3$}{k≥3}?}
Write $\sum_{x\pmod{p}}\nu(x)^k=\pi(A)N f_k$ for some function $f_k=f_k(A,M,N)$. $f_1$ is identically $1$; we will soon show that $f_2$ is $1+o(1)$. Then the upper bound \eqref{eq:HolderBound} becomes $$(\pi(A)N)^{a_1}(\pi(A)Nf_2)^{a_2}\ldots (\pi(A)Nf_m)^{a_m}Q^{1/2r},$$ say. This equals $(\pi(A)N)^{1-1/2r}f_2^{a_2}\ldots f_m^{a_m}Q^{1/2r}$.  We want to try to use higher moments of $\nu(x)$ to make an improvement, so we vary the values of $a_i$ in order to minimise $f_2^{a_2}\ldots f_m^{a_m}$. An improvement only occurs if, with respect to $A$, this is asymptotically less than $f_2^{1/2r}$; that is, if $f_3^{a_3}\ldots f_m^{a_m}$ is asymptotically smaller than $f_2^{1/2r-a_2}$. We show that this is impossible: if any of $a_3,\ldots, a_m$ are positive, then $f_3^{a_3}\ldots f_m^{a_m}$ tends to infinity while $f_2^{1/2r-a_2}$ is $O(1)$. Since $\nu(x)$ is nonnegative, we have $f_3\leq \ldots \leq f_m$, so it is enough to show that as $A\rightarrow \infty$, we have $f_3\rightarrow \infty$ while $f_2=O(1)$.

We remark here that the choice of $A$ made in Burgess's theorem ensures that $A/N=\varepsilon$ for some small positive $\varepsilon$. We may assume that $\varepsilon<1/2$.

With $c$ and $c'$ as in Lemma \ref{prop:2ndMomentUpperBound}, our second claim follows easily; upon dividing by $\pi(A)N$, the first term is $1$, the second term tends to $0$ (since $\frac{i}{p_i}\sim\frac{1}{\log{i}}$, so $\sum_{i\leq \pi(A)}\frac{i}{p_i} \sim \li{\pi(A)}\sim \frac{\pi(A)}{\log{\pi(A)}}$), and the third term (which is positive) is $$
    \varepsilon \cdot \frac{12+c'\left(\pi(A)^2+\pi(A)-16\right)}{A\pi(A)},$$ which also tends to $0$.

For the first claim, put $k=3$ into Lemma \ref{prop:exactKthMoment}: since $(P+1)^3-2P^3+(P-1)^3=6P$, $$\sum_{x\pmod{p}} \nu(x)^3=\pi(A)N+\sum_{a_1\leq A}\sum_{ca_1<a_2<a_1}\left(\left\lceil \frac{M+N}{a_1} \right \rceil - \left\lceil \frac{M}{a_2} \right\rceil \right)\cdot 6(\pi(a_1)-\pi(a_2)),$$ which we split into three parts: the diagonal term $\pi(A)N$, the off-diagonal term's main component $$6\sum_{a_1,a_2}\left(\frac{M+N}{a_1}-\frac{M}{a_2}\right)(\pi(a_1)-\pi(a_2)),$$ and its error component, whose magnitude is bounded above by $6\sum_{a_1,a_2}(\pi(a_1)-\pi(a_2)).$ In what follows, we make frequent use of the fact that if $f(n)\sim g(n)$ and $g(n)\rightarrow \infty$ as $n\rightarrow \infty$, then $\sum_{n=1}^{N}f(n)\sim \sum_{n=1}^{N}g(n)$ as $N\rightarrow \infty$. We also regularly employ Lemma \ref{prop:asympPunchConj}, which gives asymptotic equivalences for $\pi(cp_i)$ and $i-\pi(cp_i)$. All asymptotic equivalences are with respect to $A\rightarrow\infty$ unless specified.

The error component is no larger than \begin{align*}
    & 6\sum_{i=2}^{\pi(A)}\sum_{j=\pi(cp_i)+1}^{i-1}(i-j) \\
    &=6\sum_{i=2}^{\pi(A)}\left[i(i-\pi(cp_i)-1) -\left(\frac{i(i-1)}{2} - \frac{\pi(cp_i)(\pi(cp_i)+1)}{2}\right)\right] \\
    & \sim 6 \sum_{i=2}^{\pi(A)}\left[i^2(1-c) - \left(\frac{i^2}{2}-\frac{(ci)^2}{2}\right) \right] \\
    & \sim 6\left(1-c-\frac{1}{2}+\frac{c^2}{2}\right)\sum_{i=2}^{\pi(A)}i^2 \sim (c-1)^2\pi(A)^3.
\end{align*}

We now expand the summand in the main component and split it into four parts. First, we have \begin{align*}
    6\sum_{a_1,a_2} (M+N)\frac{\pi(a_1)}{a_1} &= 6(M+N)\sum_{i=2}^{\pi(A)}\sum_{j=\pi(cp_i)+1}^{i-1}\frac{i}{p_i}\\
    &= 6(M+N)\sum_{i=2}^{\pi(A)}(i-\pi(cp_i)-1)\frac{i}{p_i} \\
    &\sim 6(M+N)\sum_{i=2}^{\pi(A)}\frac{(1-c)i}{\log{i}} \\
    &\sim 6N \int_2^{\pi(A)} \frac{t}{\log{t}}\, dt= \frac{3N\pi(A)^2}{\log{\pi(A)}}.
\end{align*}

Second, we have \begin{align*}
    6\sum_{a_1,a_2} -\frac{(M+N)\pi(a_2)}{a_1}&= -6(M+N)\sum_{i=2}^{\pi(A)}\frac{1}{p_i}\sum_{j=\pi(cp_i)+1}^{i-1}j \\
    &\sim -6(M+N)\sum_{i=2}^{\pi(A)}\frac{1}{i\log{i}}\left[\frac{(i-1)i}{2}-\frac{\pi(cp_i)(\pi(cp_i)+1)}{2}\right] \\
    &\sim -3(M+N)\sum_{i=2}^{\pi(A)}\frac{i(1-c^2)}{\log{i}} \\
    &\sim -3(M+N)(1-c^2)\frac{\pi(A)^2}{2\log{\pi(A)}} \\
    &= -\frac{3N}{2}(1+c)\frac{\pi(A)^2}{\log{\pi(A)}}.
\end{align*}

Third, we have \begin{align*}
    6\sum_{a_1,a_2} -\frac{M}{a_2}\pi(a_1)&=-6M \sum_{i=2}^{\pi(A)}i \sum_{cp_i<a_2<p_i}\frac{1}{a_2} \\
    &\sim -6M\sum_{i=2}^{\pi(A)}i\left[\log{\log{p_i}}-\log{\log{cp_i}}\right].
\end{align*} The bracketed term in the summand is $$
\log{\left(\frac{\log{p_i}}{\log{c}+\log{p_i}}\right)}=-\log{\left(1+\frac{\log{c}}{\log{p_i}}\right)}\sim
 -\frac{\log{c}}{\log{p_i}}\sim -\frac{\log{c}}{\log{i}} \text{ as } i\sim\infty.$$
The sum is then asymptotic to \begin{align*}
    6M(\log{c)} \sum_{i=2}^{\pi(A)} \frac{i}{\log{i}}\sim 3M(\log{c}) \frac{\pi(A)^2}{\log{\pi(A)}}.
\end{align*}
Fourth, we have \begin{align*}
    6\sum_{a_1,a_2}\frac{M\pi(a_2)}{a_2}&=6M\sum_{i=2}^{\pi(A)}\sum_{j=\pi(cp_i)+1}^{i-1}\frac{j}{p_j}
\sim 6M \sum_{i=2}^{\pi(A)} \int_{\pi(cp_i)+1}^{i-1}\frac{1}{\log{t}}\, dt \\
&\sim 6M \sum_{i=2}^{\pi(A)}\left[\frac{i}{\log{i}}-\frac{ci}{\log{ci}}\right]\sim 6M(1-c)\sum_{i=2}^{\pi(A)}\frac{i}{\log{i}} \\
&\sim 3M(1-c)\frac{\pi(A)^2}{\log{\pi(A)}}.
\end{align*}

Combining our results, we have an error term whose magnitude is asymptotically no greater than $(c-1)^2\frac{A^3}{\log^3{A}}$, while the main term is asymptotic to $$3N\frac{A^2}{\log^3{A}}-\frac{3N}{2}(1+c)\frac{A^2}{\log^3{A}}+3M(\log{c})\frac{A^2}{\log^3{A}}+3M(1-c)\frac{A^2}{\log^3{A}}.$$ Factoring out $\pi(A)N\sim N\frac{A}{\log{A}}$ gives $f_3$. We have the diagonal term (identically $1$) and the main off-diagonal term which is asymptotic to $$\frac{3A}{\log^2{A}}\left[1-\frac{1}{2}(1+c)+\frac{M}{N}\log{c}+\frac{M}{N}(1-c)\right];$$ the term in brackets equals $$1-\frac{1}{2}-\frac{c}{2}+\frac{c}{1-c}\log{c}+c=\frac{1}{2}+\frac{c}{2}+\frac{c}{1-c}\log{c}.$$ 
We then have the error term, whose magnitude is asymptotically no greater than $(c-1)^2 \frac{A^2}{N\log^2{A}}<\varepsilon (c-1)^2 \frac{A}{\log^2{A}}$. 

The main term and error term, when combined, must be asymptotically at least $$\frac{3A}{\log^2{A}}\left[\frac{1}{2}+\frac{c}{2}+\frac{c}{1-c}\log{c}-\frac{1}{3}\varepsilon (c-1)^2\right].$$ If we can show that the coefficient of $\frac{3A}{\log^2{A}}$ is positive for every $c\in(0,1)$, then we have $f_3\rightarrow \infty$ as needed. Since $\varepsilon<1/2$, the coefficient is at least $$w(c):=\frac{1}{2}+\frac{c}{2}+\frac{c}{1-c}\log{c}-\frac{1}{6}(c-1)^2.$$ Clearly $w(c)$ is continuous for $c\in(0,1)$, and one can show that $\lim_{c\rightarrow 0^{+}} w(c)=\frac{1}{3}$ and $\lim_{c\rightarrow 1^{-}} w(c)=0$. It is therefore enough to prove that $w'(c)$ is negative for $c\in(0,1)$.

Now, $$w'(c)=\frac{11-18c+9c^2-2c^3+6\log{c}}{6(c-1)^2},$$ which has the same sign as $$11-18c+9c^2-2c^3+6\log{c}=6(1-c)+3(1-c)^2+2(1-c)^3+6\log{c}.$$ Since $\lvert c-1\rvert <1$, $\log{c}=-(1-c)-\frac{(1-c)^2}{2}-\frac{(1-c)^3}{3}-\ldots$. Thus $$6(1-c)+3(1-c)^2+2(1-c)^3+6\log{c}=6\left[-\frac{(1-c)^4}{4}-\frac{(1-c)^5}{5}-\ldots\right],$$ which is negative since $1-c>0$. This completes the proof.

\section{Burgess, \`a la Montgomery--Vaughan} \label{sec:BurgessMV}
In this section, we prove a Burgess-strength bound in the style of Montgomery--Vaughan (see \cite[Theorem 9.27]{MV2007}). However, we substitute some details from \cite{Booker2006} and make other changes. At several points in the proof, we need to rely on size assumptions about variables that will not yet have been given a value; we handle this by stating hypotheses about the variables as required, then verifying them once their values are assigned.

We need the functions defined in the following lemma in order to state our theorem; the proof will be given in Section \ref{sec:lemmas}.

\begin{table}
    \centering
    \begin{tabular}{|cc|cc|cc|cc|}
    \hline
        $c$ & $J$ & $c$ & $J$ & $c$ & $J$ & $c$ & $J$\\
    \hline
        1.891 & 3 & 0.753 & 29 & 0.338 & 149 & 0.0563 & 1447 \\
        1.388 & 5 & 0.610 & 41 & 0.254 & 227 & 0.0277 & 3299 \\
        1.036 & 11 & 0.518 & 59 & 0.203 & 347 & 0.0164 & 3457\\
        0.855 & 17 & 0.425 & 97 & 0.125 & 599 & 0 & 5393\\
    \hline
    \end{tabular}
    \caption{Values $(c,J)$ for Lemma \ref{prop:piAestimates}}
    \label{tab:AoverPiA}
\end{table}

\begin{lemma} \label{prop:piAestimates}
    Suppose $A\geq 3$ is a real number. Let $r_4(A)$ be the $c$ in Table \ref{tab:AoverPiA} corresponding to the greatest $J$ such that $J\leq A$. Then let
    \begin{align*}
    r_1(A)&=\min{\left(\frac{5}{4},\frac{\log^2{A}}{A}+\frac{20}{A^2}\right)},\, r_2(A)=\min{\left(\frac{1}{4},\frac{\log^2{A}}{A^2}+\frac{3.7}{A^2}\right)}, \\
    r_3(A)&=\min{\left(\frac{1}{2},\frac{\log{A}}{A}+\frac{4.5}{A^2}\right)}\text{ and }r_5(A)=\frac{2.55}{\log{A}}.
    \end{align*} 
    Note that each $r_i(A)$ is a decreasing function of $A$ and is $o(1)$ as $A\rightarrow \infty$.
    
    Then we have \begin{align}
        \label{prop:4.1&2&3}\frac{A}{\pi(A)^2}&\leq r_1(A),\frac{1}{\pi(A)^2}\leq r_2(A), \frac{1}{\pi(A)}\leq r_3(A),\\
        \label{prop:4.4&5}\frac{A}{\pi(A)}-\log{A}+1&\leq r_4(A), \text{ and } -1+\frac{A}{\pi(A)^2}\sum_{i\leq \pi(A)}\frac{i}{p_i}\leq r_5(A).   
    \end{align}
\end{lemma}
\begin{theorem} \label{prop:BurgessMV}
    Let $p$ be a prime number greater than $10^3$. Let $M_0$ and $N_0$ be positive integers with $N_0\leq p^{1/2}$ and select some integer $r\in[2,15]$. Fix some real $k\in (0,2^{-r/(r-1)})$. Write $c_0=\frac{M_0}{M_0+N_0}$ and let $c'_0=\max(1-c_0,13/61)$. Let $$
B_0= \left[\frac{1}{2}(2r-1)!!\right]^{1/r}p^{1/2r},
A_0=\frac{kp^{(r+1)/4r}}{\lceil B_0 \rceil}, \text{ and }
\varepsilon= \frac{k}{\lceil B_0 \rceil}.$$ Let $C_0=\frac{p^{\frac{1}{2}+\frac{1}{2r}}}{k}\left[\frac{1}{2}(2r-1)!!\right]^\frac{1}{r}$ and 
$$E= r_4(A_0)+\frac{244}{105}r_1(A_0)+2c'_0 r_5(A_0)+12\varepsilon r_2(A_0)+\varepsilon c'_0 r_3(A_0),$$ where the $r_i$ are as in Lemma \ref{prop:piAestimates}, and let 
$$Y=\left(\log{\frac{kN_0}{\left(\frac{1}{2}(2r-1)!!\right)^{1/r}}-\frac{\log{p}}{2r}}-1+2c_0'+\varepsilon c_0'+E\right)\left[2r\left(\frac{1}{2}(2r-1)!!\right)^{1/r}\right].$$

Suppose that $M_0+\frac{N_0}{1-k}\leq C_0$, $A_0\geq 3$ and $\lceil B_0 \rceil \geq 2r$. 
Then $$\lvert S(M_0,N_0)\rvert \leq \frac{\max(1,Y/k)^{1/2r}}{1-2k^{1-1/r}}p^\frac{r+1}{4r^2}N_0^{1-1/r}.$$
\end{theorem}
\begin{rem}
    We will see later that the choice $k=(4r-2)^{\frac{-r}{r-1}}$ is nearly optimal.
\end{rem}
\begin{proof}
Clearly, our bound is at least $p^\frac{r+1}{4r^2}N_0^{1-1/r}$, so it is implied by the trivial bound $N_0$ if $N_0\leq p^\frac{r+1}{4r}$. In the ensuing proof, we allow $M$ and $N$ to take on any positive integer values satisfying the following conditions: $p^{\frac{r+1}{4r}}< N\leq N_0$, $\frac{M}{M+N}\geq c_0$, $M+\frac{N}{1-k}\leq C_0$. We also write $c=\frac{M}{M+N}$ and $c'=\max(1-c,13/61)$.

First, we shift the character sum by a positive integer $ab$.
\begin{hyp} \label{hyp:hypothesis1}
    Suppose that $ab\leq kN$.
\end{hyp}
Then \begin{align*}
    S(M,N)&=\sum_{n\in[M,M+N)}\chi(n+ab)+\sum_{n\in[M,M+ab)}\chi(n)-\sum_{n\in[M+N,M+N+ab)}\chi(n)\\
    &=\sum_{n\in[M,M+N)}\chi(n+ab)+S(M,ab)-S(M+N,ab).
\end{align*}
Now define $$J(y)=\max_{\substack{N\leq y\\ M>0\\ \frac{M}{M+N}\geq c_0 \\ M+\frac{N}{1-k}\leq C_0}}\lvert S(M,N)\rvert.$$ $J$ is an increasing function of $y$, because once some $(M,N)$ satisfies the constraints for a certain $y$, it will certainly satisfy them for any $y'\geq y$. Observe that $$\frac{M}{M+ab}>\frac{M}{M+N}\geq c_0 \text{ and }\frac{M+N}{M+N+ab}> \frac{M+N}{M+2N}>\frac{M}{M+N}\geq c_0.$$ Further, we have $$M+\frac{ab}{1-k}<M+N+\frac{ab}{1-k} \leq M+N+\frac{kN}{1-k}=M+\frac{N}{1-k}\leq C_0.$$ It follows that both $\lvert S(M,ab)\rvert $ and $\lvert S(M+N,ab)\rvert $ are bounded above by $J(ab)$. Thus $$S(M,N)= \sum_{n\in[M,M+N)}\chi(n+ab)+2\theta J(ab),$$ for some $\lvert \theta \rvert\leq 1$. We sum this equation over primes $a\leq A$ and integers $1\leq b \leq B$: $$\pi(A)B\cdot S(M,N) = \sum_{\substack{1\leq a \leq A\\ a \text{ prime}}}\sum_{1\leq b \leq B}\sum_{n\in [M,M+N)} \chi(n+ab)+2\theta' \pi(A)BJ(AB),$$ where similarly $\lvert \theta' \rvert \leq 1$. 
\begin{hyp} \label{hyp:hypothesis2}
    Suppose that $A<p$.
\end{hyp}
Let $\nu(x)$ be the number of representations of $x$ as $a^{-1}n\pmod{p}$ where $a$ is a prime in $[1,A]$ and $n$ is an integer in $[M,M+N)$. Abbreviating the triple sum above as $\sum_{n,a,b}$, we have \begin{align*}
    \left\lvert \sum_{n,a,b}\chi(n+ab) \right\rvert &=\left \lvert \sum_{x=1}^{p} \sum_{\substack{n,a\\ n\equiv{ax\pmod{p}}}}\chi(a)\sum_b{\chi(x+b})\right\rvert \\
    &\leq \sum_{x=1}^{p} \nu(x)\left\lvert \sum_b{\chi(x+b)}\right\rvert \\
    &\leq \left(\sum_{x=1}^{p}\nu(x)^{2r/(2r-1)}\right)^{1-1/2r}\left(\sum_{x=1}^{p}\left\lvert \sum_b \chi(x+b)\right\rvert^{2r} \right)^{1/2r}\\
    &\leq \left(\sum_{x=1}^{p} \nu(x)\right)^{1-1/r}\left(\sum_{x=1}^{p}\nu(x)^2\right)^{1/2r}\left(\sum_{x=1}^{p}\left\lvert \sum_b \chi(x+b)\right\rvert^{2r} \right)^{1/2r}.
\end{align*} The final two steps each involve an application of H\"older's inequality.

Now, $\sum_{x=1}^p \nu(x)=\pi(A)N$. Assume for now that we know some $f=f(A,M,N)$ such that $\sum_{x=1}^{p} \nu(x)^2\leq \pi(A)Nf$; later we will use Lemma \ref{prop:2ndMomentUpperBound} to find such an $f$. The third sum can be dealt with as in \cite{Booker2006}; we get $$
    \sum_{x=1}^{p}\left\lvert \sum_{b=1}^{B} \chi(x+b)\right\rvert^{2r} = \sum_{1\leq b_1,\ldots,b_{2r}\leq B}\sum_{x=1}^{p} \chi((x+b_1)\ldots (x+b_{2r}))$$ because $\chi$ has order $2$. Examining the polynomial $(x+b_1)\ldots(x+b_{2r})$, we consider three cases.
    \begin{hyp} \label{hyp:hypothesis3}
        Suppose that $B\geq 2r$.
    \end{hyp}
    \begin{enumerate}
        \item Suppose that the polynomial is a square. In this case, the sum over $x$ is at most $p$. In how many ways can this happen? Choose $b_1$; there are $B$ ways to do this. Then, choose some $2\leq \ell \leq 2r$ and set $b_\ell=b_1$; there are $2r-1$ ways to do this. Move to the first $b_i$ that is not chosen; choose it and its twin, and continue until all the $b_i$ are chosen. This gives $(2r-1)!!\cdot B^r$ ways.
        \item Suppose that all of the $b_i$ are distinct; this can occur in exactly $\frac{B!}{(B-2r)!}$ ways. By Lemma \ref{prop:Weil}, the sum over $x$ has magnitude no greater than $2(r-1)\sqrt{p}+1$. 
        \item Suppose that the polynomial is not a square, but not all of the $b_i$ are distinct. This can occur in no more than $B^{2r}-\frac{B!}{(B-2r)!}$ ways. Since $\chi$ has order $2$, the polynomial's image under the character will not change if we delete a repeated factor $(x+c)^2$, except possibly at $-c$. The new polynomial has degree $2(r-1)$, so Lemma \ref{prop:Weil} (with an adjustment at $x=-c$) now gives an upper bound of $ 2(r-2)\sqrt{p}+2$.
    \end{enumerate}
    
    This gives a total contribution no greater than \begin{align*}&\left[2(r-2)\sqrt{p}+2\right]\left(B^{2r}-\frac{B!}{(B-2r)!}\right)+[2(r-1)\sqrt{p}+1]\left(\frac{B!}{(B-2r)!}\right)+(2r-1)!!\cdot p B^r \\
    &=2(r-1)B^{2r}\sqrt{p}+(2r-1)!!\cdot  B^r p-\left[(2\sqrt{p}-1)\left(B^{2r}-\frac{B!}{(B-2r)!}\right)-B^{2r}\right],
    \end{align*} which we will denote by $W(B,r,p)$, or just $W$ for short. Now, let $AB=kN$; this verifies Hypothesis \ref{hyp:hypothesis1}. Then $\pi(A)=\frac{kN\pi(A)}{AB}$. Combining our results, we find that \begin{align*}
    \lvert S(M,N) \rvert &\leq \frac{1}{\pi(A)B}\left\lvert \sum_{n,a,b}\chi(n+ab) \right\rvert +2J(AB) \\
    &\leq \frac{(\pi(A)N)^{1-1/2r}f^{1/2r}W^{1/2r}}{\pi(A)B} +2J(AB) \\
    &=\left(\frac{Af}{\pi(A)}\right)^{1/2r}k^{-1/2r}N^{1-1/r}\left(\frac{W}{B^{2r-1}}\right)^{1/2r}+2J(AB).
\end{align*} For now, ignore the last term of $W$, so that $$\frac{W(B,r,p)}{B^{2r-1}}\approx 2(r-1)B\sqrt{p}+(2r-1)!!\cdot B^{1-r}p=:CB+DB^{1-r},$$ where $C$ and $D$ are constant with respect to $B$. The function $F(X):=CX+DX^{1-r}$ is minimised at $$(D(r-1)/C)^{1/r}=\left((2r-1)!!/2\right)^{1/r}p^{1/2r}=B_0.$$ Choose $B=\ceil{B_0}$, so that $B-B_0<1$ and Hypothesis \ref{hyp:hypothesis3} is validated. Now, $$F(B_0)=C^{1-\frac{1}{r}}D^{\frac{1}{r}}r(r-1)^{\frac{1}{r}-1}=2r\left[\frac{1}{2}(2r-1)!!\right]^{1/r}p^{\frac{1}{2}+\frac{1}{2r}};$$ we remark here that the coefficient $2r$ differs from the $4(r-1)$ stated at \cite[(34)]{Booker2006}. The author believes that $2r$ is correct.

By how much can $F(B)$ exceed this? We estimate the difference by using Taylor's theorem in the following form:
$$F(B)-F(B_0)=F'(B_0)(B-B_0)+\frac{F''(\xi)}{2!}(B-B_0)^2,$$ for some $\xi$ between $B_0$ and $B$. $F'(B_0)=0$ and $F''(\xi)=Dr(r-1)\xi^{-r-1}$, so \begin{align}F(B)-F(B_0)&< \frac{1}{2}F''(\xi)\leq \frac{p}{2}(2r-1)!!\cdot r(r-1)B_0^{-r-1}.
\end{align}

Actually, including the last term of $W$ will be enough to cancel this out; once divided by $B^{2r-1}$, the term from $W$ is $$\frac{2\sqrt{p}-1}{B^{2r-1}}\left[B^{2r}-\frac{B!}{(B-2r)!}\right]-B=(2\sqrt{p}-1)B\left[1-\frac{B!}{B^{2r}(B-2r)!}\right]-B.$$ The bracketed part is \begin{align*}
    &1-\frac{(B)(B-2r+1)(B-1)(B-2r+2)\ldots (B-r+1)(B-r)}{B^{2r}}\\ &\geq 1-\frac{(B-r+1/2)^{2r}}{B^{2r}}= 1-\left(1-\frac{r-1/2}{B}\right)^{2r} \\
    &=1-\exp{\left(2r \log{\left(1-\frac{r-1/2}{B}\right)}\right)}\geq 1-\exp{\left(\frac{-r(2r-1)}{B}\right)},
\end{align*} because $\log{(1-x)}\leq -x$. It therefore suffices to show that $$B_0\left[(2\sqrt{p}-1)\left(1-\exp\left(\frac{-r(2r-1)}{B_0+1}\right)\right)-1\right]\geq \frac{p}{2}(2r-1)!!\cdot r(r-1)B_0^{-r-1}$$ or, equivalently, that \begin{align}\label{eq:errorTermSizes1}B_0^2 \left[1-\exp\left(\frac{-r(2r-1)}{B_0+1}\right)\right]\geq \frac{B_0^2}{2\sqrt{p}-1}+\frac{r(r-1)}{2-\frac{1}{\sqrt{p}}}.\end{align} Since $p> 10^3$, the right-hand side is less than $$\frac{\left[\frac{1}{2}(2r-1)!!\right]^{2/r}\left(10^3\right)^{1/r-1/2}}{2-\frac{1}{\sqrt{10^3}}}+\frac{r(r-1)}{2-\frac{1}{\sqrt{10^3}}},$$ so it is enough to prove that the left-hand side exceeds this. The derivative of the left-hand side of \eqref{eq:errorTermSizes1} with respect to $B_0$ is $$\left(2-2e^{\frac{-r(2r-1)}{B_0+1}}\right)B_0-e^\frac{-r(2r-1)}{B_0+1}r(2r-1)\frac{B_0^2}{(B_0+1)^2}.$$ If we write $\kappa:=\frac{r(2r-1)}{B_0+1}>0$, then this is $$(2-2e^{-\kappa})B_0-\kappa e^{-\kappa}\cdot \frac{B_0^2}{B_0+1}\geq 2B_0-2e^{-\kappa}B_0-\kappa e^{-\kappa}B_0=B_0[2-2e^{-\kappa}-\kappa e^{-\kappa}].$$ The Maclaurin expansion of $e^\kappa$ implies that $e^\kappa> 1+\kappa/2$ for $\kappa>0$; it follows that $2-2e^{-\kappa}-\kappa e^{-\kappa}> 0$ and so the left-hand side of \eqref{eq:errorTermSizes1} is increasing with respect to $B_0$, which is in turn increasing with respect to $p$. Thus, if \eqref{eq:errorTermSizes1} holds at $p=p_0$ for some $p_0$, it will hold for $p\geq p_0$. It is readily verified, then, that \eqref{eq:errorTermSizes1} holds whenever $p\geq 10^3$ and $2\leq r \leq 15$.

Substituting $\frac{W(B,r,p)}{B^{2r-1}}\leq F(B_0)$ gives $$\lvert S(M,N)\rvert \leq \left(\frac{Af}{\pi(A)}\right)^{1/2r}k^{-1/2r}N^{1-1/r}p^{\frac{r+1}{4r^2}}\left[2r\left(\frac{1}{2}(2r-1)!!\right)^{1/r}\right]^{1/2r}+2J(AB).$$

Next, we examine the $\frac{Af}{\pi(A)}$ term. We have $A=\frac{kN}{B}>A_0\geq 3$ by the definitions of $B$ and $A_0$ and our assumptions about $N$ and $A_0$. Similarly, we have $\frac{A}{N}=\frac{k}{B}=\varepsilon$.
\begin{hyp} \label{hyp:hypothesis4}
    Suppose that $A(M+N)\leq p$.
\end{hyp}
By Lemma \ref{prop:2ndMomentUpperBound}, \begin{align*} f&\leq 1+\frac{2}{\pi(A)}\left[\frac{122}{105}+c'\sum_{5\leq i \leq \pi(A)}\frac{i}{p_i}\right]+\frac{12+c'[\pi(A)^2+\pi(A)-16]}{\pi(A)N} \\
\frac{Af}{\pi(A)}&\leq \frac{A}{\pi(A)}+\frac{2A}{\pi(A)^2}\left[\frac{122}{105}+c'\sum_{5\leq i \leq \pi(A)}\frac{i}{p_i}\right]+\frac{A(12+c'[\pi(A)^2+\pi(A)-16])}{\pi(A)^2N}.
\end{align*}
We may now write \begin{align}\label{eq:AfoverPiA}\frac{Af}{\pi(A)}= \log{A}-1+2c'+\varepsilon c'+E(A,c',\varepsilon),\end{align} where $E(A,c',\varepsilon)$ equals \begin{align*}&\left[\frac{A}{\pi(A)}-\log{A}+1\right]+\frac{\frac{244}{105}A}{\pi(A)^2}+2c'\left[-1+\frac{A}{\pi(A)^2}\sum_{5\leq i \leq \pi(A)}\frac{i}{p_i}\right]+\frac{\varepsilon(12-16 c')}{\pi(A)^2}+\frac{\varepsilon c'}{\pi(A)}\\
&< r_4(A)+\frac{244}{105}r_1(A)+2c'r_5(A)+12\varepsilon r_2(A)+\varepsilon c' r_3(A),
\end{align*} by Lemma \ref{prop:piAestimates}. Since each $r_i$ is decreasing with $A$ and we have $A>A_0$, $c'\leq c_0'$, this is no greater than $E$.

Also, $$\log{A}=\log{\frac{kN}{B}} 
\leq\log{\frac{kN_0}{\left(\frac{1}{2}(2r-1)!!\right)^{1/r}}}-\frac{\log{p}}{2r}.$$

Using $c'\leq c_0'$ again, we get \begin{align}\notag  \lvert S(M,N)\rvert &\leq \left(\log{\frac{kN_0}{\left(\frac{1}{2}(2r-1)!!\right)^{1/r}}-\frac{\log{p}}{2r}}-1+2c_0'+\varepsilon c_0'+E\right)^{1/2r} \\ \notag &\times k^{-1/2r}N^{1-1/r}p^{\frac{r+1}{4r^2}}\left[2r\left(\frac{1}{2}(2r-1)!!\right)^{1/r}\right]^{1/2r}+2J(kN)\\
\notag &=Y^{1/2r}k^{-1/2r}N^{1-1/r}p^\frac{r+1}{4r^2}+2J(kN) \\
\label{eq:altogether} &\leq \max(1,Y/k)^{1/2r}N^{1-1/r}p^\frac{r+1}{4r^2}+2J(kN).
\end{align} 

We now put an upper bound on $J(N)$. If $N\leq p^\frac{r+1}{4r}$, then $$N^{1-1/r}p^\frac{r+1}{4r^2}\geq N\cdot p^{-\frac{r+1}{4r^2}}\cdot p^{\frac{r+1}{4r^2}}=N.$$ Combining this with the trivial bound, we see that $J(N)$ cannot exceed the right-hand side of \eqref{eq:altogether}.

Suppose instead that $N\in (p^\frac{r+1}{4r}, N_0]$, and examine only those $(M_1,N_1)$ considered in $J(N)$. If $N_1\in (p^\frac{r+1}{4r},N_0]$, then we can apply \eqref{eq:altogether}. If $N_1\leq p^{\frac{r+1}{4r}}$, then $S(M_1,N_1)\leq p^{\frac{r+1}{4r}}$. But we also have $N^{1-1/r}p^\frac{r+1}{4r^2}>p^{\frac{r+1}{4r}-\frac{r+1}{4r^2}+\frac{r+1}{4r^2}}=p^\frac{r+1}{4r}$. In both cases, $S(M_1,N_1)$ cannot exceed the right-hand side of \eqref{eq:altogether}.

Therefore for any $N\leq N_0$, \begin{align}\label{eq:Jc(N)}J(N)\leq Zk^{-1/2r}N^{1-1/r}+2J(kN),\end{align} where we have written $Z:=p^\frac{r+1}{4r^2}\max(k,Y)^{1/2r}$.

Applying \eqref{eq:Jc(N)} repeatedly $Q$ times yields
\begin{align*}J(N)&\leq Zk^{-1/2r}N^{1-1/r}+\ldots+2^QZk^{-1/2r}(k^Q N)^{1-1/r}+2^{Q+1}J(k^{Q+1}N) \\
&=Zk^{-1/2r}N^{1-1/r}\left[1+2k^{1-1/r}+\ldots+2^Q (k^Q)^{1-1/r}\right]+2^{Q+1}J(k^{Q+1}N)\\
&\leq ZN^{1-1/r}\frac{k^{-1/2r}}{1-2k^{1-1/r}}+2^{Q+1}k^{Q+1}N,
\end{align*} where we have used $2k^{1-1/r}<1$ and the trivial bound for $J$. The function $\frac{x^{-1/2r}}{1-2x^{1-1/r}}$ has a minimum at $x=(4r-2)^{-\frac{r}{r-1}}$ at which it equals $\frac{(4r-2)^\frac{2r-1}{2r-2}}{4(r-1)}.$ A sensible choice of $k$ is therefore $(4r-2)^\frac{-r}{r-1}$ (which is, as required, less than $2^{-\frac{r}{r-1}}$). We note, though, that because $Z$ is a function of $k$, this may not be optimal. In any case, since $k<1/2$, $2^{Q+1}k^{Q+1}N$ can be made arbitrarily small by choosing $Q$ to be large. Thus $J(N)\leq \frac{k^{-1/2r}}{1-2k^{1-1/r}}ZN^{1-1/r};$ taking $N=N_0$ yields the theorem.

It remains to check those hypotheses that have not yet been validated. We begin with Hypothesis \ref{hyp:hypothesis2} (that $A<p$); we have $$A=\frac{kN}{B}< \frac{2^{-r/(r-1)}}{\left(\frac{1}{2}(2r-1)!!\right)^{\frac{1}{r}}}\times \frac{N_0}{p^{\frac{1}{2r}}}.$$ We know that $N_0/p^{1/2r}<p$; for $2\leq r \leq 15$ it is readily verified that $\frac{2^{-r/(r-1)}}{\left(\frac{1}{2}(2r-1)!!\right)^{\frac{1}{r}}}<1.$

For Hypothesis \ref{hyp:hypothesis4}, it is enough to show that $M+N\leq \frac{p}{A}=\frac{pB}{kN}$. First, recall that $M+N<M+\frac{N}{1-k}\leq C_0$. Next, we have $$\frac{pB}{kN}\geq \frac{p\left[\frac{1}{2}(2r-1)!!\right]^{1/r}p^{1/2r}}{kp^{1/2}}=\frac{p^{\frac{1}{2}+\frac{1}{2r}}}{k}\left[\frac{1}{2}(2r-1)!!\right]^{1/r}=C_0.$$

\end{proof}

\section{Burgess, \`a la Booker}\label{sec:BurgessBooker}

We now prove Theorem \ref{prop:BurgessBooker}; this time, we follow \cite{Booker2006} very closely. There are slight differences, including:
\begin{enumerate}
    \item we restrict the set of intervals that can be reached during the shifting process,
    \item we employ Lemma \ref{prop:2ndMomentUpperBound} to get a better upper bound on the second moment of $\nu(x)$, and
    \item we correct an error in \cite{Booker2006}.
\end{enumerate}
\begin{proof}[Proof of Theorem \ref{prop:BurgessBooker}]
We use strong induction: suppose we know that $\lvert S(M,N)\rvert \leq K p^\frac{r+1}{4r^2}N^{1-1/r}$ for every pair of positive integers $(M,N)$ satisfying the following: 
\begin{align}\label{eq:conditions} (M,N)\text{ is such that }N< N_0, \frac{M}{M+N}\geq c_0 \text{ and } M+\frac{N}{1-k}\leq C_0.\end{align}
    
     It is enough to prove that the same holds if we relax $N<N_0$ to $N\leq N_0$. Since $K>1$, the trivial bound of $N_0$ implies this if $N_0\leq p^\frac{r+1}{4r}$, so assume that $N_0>p^\frac{r+1}{4r}$.
    
    We will shift by an integer $ab$.
    \begin{hyp} \label{hyp:hypothesis5}
        Suppose that $\lvert ab \rvert \leq kN_0$.
    \end{hyp}
    Note that $k<1$. If $ab\geq 0$, then $$S(M_0,N_0)=\sum_{n\in[M_0,M_0+N_0)}\chi(n+ab)+S(M_0,ab)-S(M_0+N_0,ab)$$ as before. Otherwise, $ab<0$, and
    $$S(M_0,N_0)=\sum_{n\in[M_0,M_0+N_0)}\chi(n+ab)-S(M_0-\lvert ab\lvert ,\lvert ab\rvert )+S(M_0+N_0-\lvert ab\rvert,\lvert ab\rvert).$$

    As before, $(M_0,ab)$ and $(M_0+N_0,ab)$ satisfy \eqref{eq:conditions}. We now verify the $ab<0$ case: $$\frac{M_0-\lvert ab\rvert }{M_0}>\frac{M_0}{M_0+N_0}\geq c_0$$ since $\lvert ab \rvert \leq kN_0$ and $\frac{M_0}{M_0+N_0}\geq k$. Further, $$\frac{M_0+N_0-\lvert ab\rvert }{M_0+N_0}>\frac{M_0}{M_0+N_0}\geq c_0.$$ One can show that $$M_0-\lvert ab\rvert +\frac{\lvert ab\rvert}{1-k}\text{ and }M_0+N_0-\lvert ab\rvert +\frac{\lvert ab\rvert }{1-k}$$ do not exceed $M_0+\frac{N_0}{1-k}\leq C_0$. Finally, $M_0+N_0-\lvert ab \rvert$ is clearly positive; $M_0-\lvert ab\rvert$ is too because $\lvert ab \rvert \leq kN_0 \leq (1-k)M_0<M_0$. Therefore $(M_0-\lvert ab\rvert ,\lvert ab \rvert)$ and $(M_0+N_0-\lvert ab\rvert ,\lvert ab \rvert)$ both satisfy \eqref{eq:conditions}.

    We may then use the inductive hypothesis as we average over primes $\lvert a \rvert\leq A$ and integers $\lvert b \rvert \leq B$; this gives $$\left \lvert S(M_0,N_0)-\frac{\sum_{n,a,b}\chi(n+ab)}{2\pi(A)(2B+1)}  \right \rvert \leq \frac{2K p^\frac{r+1}{4r^2}[A(B+1/2)]^{1-\frac{1}{r}}}{(2-1/r)^2}\frac{\sum_a a^{1-1/r}}{\pi(A)A^{1-1/r}/(2-1/r)}.$$

    As before, we assume that each $a$ has an inverse modulo $p$ and that there are enough options for each $b_i$ in the polynomial $(x+b_1)\ldots(x+b_{2r})$ for our counting argument to be valid.
    \begin{hyp} \label{hyp:hypothesis6}
        Suppose that $A<p$.
    \end{hyp}
    \begin{hyp} \label{hyp:hypothesis7}
        Suppose that $2B+1\geq 2r$.
    \end{hyp}
    Let $\nu(x)$ be the number of representations of $x$ as $a^{-1}n\pmod{p}$, where $\lvert a\rvert$ is a prime number in $[1,A]$ and $n$ is an integer in $[M_0,M_0+N_0)$. We have $\sum_{x=1}^{p}\nu(x)=2\pi(A)N_0$ and will write $\sum_{x=1}^{p}\nu(x)^2\leq 2\pi(A)N_0f$. Set $A(B+1/2)=kN_0$ (verifying Hypothesis \ref{hyp:hypothesis5}), so that \begin{align*}
        \frac{\sum_{n,a,b}\chi(n+ab)}{2\pi(A)(2B+1)}&\leq \frac{(2\pi(A)N_0)^{1-1/r}(2\pi(A)N_0 f)^{1/2r}W^{1/2r}}{2\pi(A)(2B+1)} \\
        &=\left(\frac{Af}{4k\pi(A)}\right)^{1/2r} {N_0}^{1-1/r}\left(\frac{W}{(2B+1)^{2r-1}}\right)^{1/2r},
    \end{align*} where we have written \begin{align*}W&=2(r-1)(2B+1)^{2r}\sqrt{p}+(2r-1)!!(2B+1)^r p-\\ &\left[(2\sqrt{p}-1)\left((2B+1)^{2r}-\frac{(2B+1)!}{(2B+1-2r)!}\right)-(2B+1)^{2r}\right].\end{align*}
    Choose $2B+1$ to be the least odd integer greater than or equal to $\left[\frac{1}{2}(2r-1)!!\right]^{1/r}p^{1/2r}$ (this validates Hypothesis \ref{hyp:hypothesis7}). As in Theorem \ref{prop:BurgessMV}, one may verify that $\frac{W}{(2B+1)^{2r-1}}\leq 2r\left[\frac{1}{2}(2r-1)!!\right]^{1/r}p^{\frac{1}{2}+\frac{1}{2r}}$. As mentioned above, a coefficient of $4(r-1)$ is given in \cite{Booker2006}; the author believes that $2r$ is instead correct.

    This yields \begin{align*}
        \lvert S(M_0,N_0)\rvert \leq \left(\frac{2r\left[\frac{1}{2}(2r-1)!!\right]^{1/r}Af}{4k\pi(A)}\right)^{1/2r}N_0^{1-1/r}p^\frac{r+1}{4r^2}+\frac{2Kp^\frac{r+1}{4r^2}(kN_0)^{1-1/r}g}{(2-1/r)^2},
    \end{align*} where $g=\frac{\sum_{a\text{ prime}\leq A} a^{1-1/r}}{\pi(A)A^{1-1/r}/(2-1/r)}$. We need the right-hand side to be at most $Kp^\frac{r+1}{4r^2}N_0^{1-1/r}$, this happens if
    \begin{align}\label{eq:Kcondition1}
        K\geq \left(\frac{\frac{r}{2}\left[\frac{1}{2}(2r-1)!!\right]^{1/r}Af}{\pi(A)}\right)^{1/2r}\frac{k^{-1/2r}}{1-\frac{2}{(2-1/r)^2}k^{1-1/r}g}.
    \end{align} Substituting $k=\left(\frac{2r-1}{2r^2}\right)^{r/(r-1)}$, we can rewrite \eqref{eq:Kcondition1} as
    \begin{align}
        \label{eq:WinIfK}K&\geq \tilde{\alpha}(r)\left[\left(\frac{2r}{r-1}\cdot \frac{Af}{\pi(A)}\right)\left(\frac{2r-2}{2r-1-g}\right)^{2r}\right]^{1/2r}, \text{ where}\\
        \notag \tilde{\alpha}(r)&=\left(\frac{r-1}{4}\right)^{1/2r}\left[\frac{1}{2}(2r-1)!!\right]^{1/2r^2}\left(\frac{2r^2}{2r-1}\right)^{1/(2r-2)}\frac{2r-1}{2r-2}.
    \end{align} We note that $\tilde{\alpha}(r)<\alpha(r)$.

    Recall that the $\nu(x)$ from Lemma \ref{prop:exactKthMoment} has $\mathcal{A}$ containing the primes up to $A$. As remarked after the proof of Lemma \ref{prop:exactKthMoment}, we may assume that $A(M_0+N_0)\leq p/2$ and double the upper bound.
    \begin{hyp}\label{hyp:hypothesis8}
        Suppose that $A(M_0+N_0)\leq p/2$.
    \end{hyp}

    We wrote $\sum_{x=1}^{p}\nu(x)^2\leq 2\pi(A)N_0f$, so our $f$ and $\frac{Af}{\pi(A)}$ will be the same as in Theorem \ref{prop:BurgessMV}. One can verify that $\frac{A}{N_0}<0.0017$ for any permitted $r$ and $p$, so by \eqref{eq:AfoverPiA}, we get $\frac{Af}{\pi(A)}\leq \log{A}-0.9983+2c_0'+E(A,c_0',\varepsilon),$ with $\varepsilon=\frac{A}{N_0}$. Assume for now that $A\geq 10^5$; along with $c_0'<1$ and $\varepsilon<0.0017$, it is not difficult to show that $E(A,c_0',\varepsilon)<\frac{5.14}{\log{A}}$.
    
    We also have $$\log{A}\leq \left(\frac{r-1}{2r}\right)\log{p}+\log{\frac{4\left(\frac{2r-1}{2r^2}\right)^{r/(r-1)}}{\left[\frac{1}{2}(2r-1)!!\right]^{1/r}}}.$$ Combining yields \begin{align*}\frac{2r}{r-1}\cdot \frac{Af}{\pi(A)}&\leq \log{p}+\frac{2r}{r-1}\left(\log{\frac{4\left(\frac{2r-1}{2r^2}\right)^{r/(r-1)}}{\left[\frac{1}{2}(2r-1)!!\right]^{1/r}}}-0.9983+2c_0'\right)+\frac{2r\cdot 5.14}{(r-1)\log{A}}\\
    &=\log{p}+\tilde{\beta}(r)+\frac{4r}{r-1}c_0'+\frac{10.28\frac{r}{r-1}}{\log{A}}=D(p,r,c_0')+\frac{E(r)}{\log{A}},
    \end{align*} say. Note that $\tilde{\beta}(r)<\beta(r)$. Taking $r=2$ as an example, \eqref{eq:WinIfK} and Lemma \ref{prop:eq47Part2} imply that it is enough to have $$K\geq \alpha(2)\left[\left(D(p,2,c_0')+\frac{E(2)}{\log{A}}\right)\left(1-\frac{2/3}{\log{A}}+\frac{1.23228}{\log^2{A}}\right)\right]^{1/2r}.$$ Now $E(2)=20.56$, and because $p>10^{23}$ and $c_0'\geq 13/61$ we have $D(p,2,c_0')> 47.55$. The expression inside the square brackets is less than \begin{align*}&D(p,2,c_0')+\frac{20.56}{\log{A}}\left(1-\frac{2/3}{\log{A}}\right)-\frac{47.55}{\log{A}}\left(\frac{2}{3}-\frac{1.23228}{\log{A}}\right)+\frac{25.34}{\log^3{A}}\\
    &< D(p,2,c_0')+\frac{1}{\log{A}}\left(-11.14+\frac{44.89}{\log{A}}+\frac{25.34}{\log^2{A}}\right)<D(p,2,c_0'),
    \end{align*} because $A\geq 10^5$. So for these $A$, it is enough to take $K\geq \alpha(2)D(p,2,c_0')^{1/4}$. A similar argument establishes that $K\geq \alpha(r)D(p,r,c_0')^{1/2r}$ is sufficient for any $2\leq r \leq 15$. We now need to check each smaller value of $A$ that may occur; we have $$A=\frac{2kN}{2B+1}>\frac{2\left(\frac{2r-1}{2r^2}\right)^{\frac{r}{r-1}}p^\frac{r+1}{4r}}{\left[\frac{1}{2}(2r-1)!!\right]^{1/r}p^{1/2r}+2}>\frac{2\left(\frac{2r-1}{2r^2}\right)^{\frac{r}{r-1}}(10^{23})^\frac{r+1}{4r}}{\left[\frac{1}{2}(2r-1)!!\right]^{1/r}(10^{23})^{1/2r}+2}=:A_r.$$

    We cover the range $(A_r,10^5)$ with the intervals $[p_{\pi(A_r)},p_{\pi(A_r)+1}),\ldots,[p_{9592},p_{9593})$ and consider each separately. When $A\in [p_k,p_{k+1})$, we have \begin{align*}
        \frac{2r}{r-1}\frac{A}{\pi(A)}f &< \frac{2r}{r-1}\cdot \frac{p_{k+1}}{k}\left(1+\frac{2\left[\frac{122}{105}+\sum_{5\leq i \leq k}\frac{i}{p_i}\right]}{k}+0.0017\cdot \frac{k^2+k-4}{kp_k}\right)\\
        \left(\frac{2r-2}{2r-1-g}\right)^{2r}&\leq \left(\frac{2r-2}{2r-1-\frac{(2-1/r)\sum_{1\leq i \leq k}p_i^{1-1/r}}{kp_k^{1-1/r}}}\right)^{2r}.
    \end{align*} For each $r$ and for all $\pi(A_r)\leq k \leq 9592$, we use these upper bounds to verify that $D(10^{23},r,13/61)\geq \left(\frac{2r}{r-1}\cdot \frac{Af}{\pi(A)}\right)\left(\frac{2r-2}{2r-1-g}\right)^{2r}$, which suffices because we must have $D(p,r,c_0')>D(10^{23},r,13/61)$.

    As for Hypotheses \ref{hyp:hypothesis6} and \ref{hyp:hypothesis8}, they can be validated in the same way as Hypotheses \ref{hyp:hypothesis2} and \ref{hyp:hypothesis4} (see the proof of Theorem \ref{prop:BurgessMV}).
      
\end{proof}

\section{Lemmas}\label{sec:lemmas}

\begin{lemma}{\cite[Corollary 5.2]{Dusart2018}}\label{prop:DusartPi}
\begin{enumerate}
    \item If $x\geq 11813$, then $\pi(x)=\frac{x}{\log{x}}+\frac{x}{\log^2{x}}+\frac{kx}{\log^3{x}}$ for some $k\in [1.46184,2.53816]$.
    \item If $x\geq 599$, then $\pi(x)=\frac{x}{\log{x}}+\frac{\ell x}{\log^2{x}}$ for some $\ell\in [1,1.2762]$.
\end{enumerate}
 
\end{lemma}
\begin{lemma}\label{prop:AoverPiA} Suppose that $A\geq 10^5$. Then $\frac{A}{\pi(A)}=\log{A}-1+\frac{\ell}{\log{A}}$, for some $\ell\in (-1.53816,-0.30273)$.\end{lemma}
\begin{proof}
    Let $k$ be as in Lemma \ref{prop:DusartPi}. Then $$\log(A)\left[\frac{A}{\pi(A)}-\log{A}+1\right]=\frac{1-k+\frac{k}{\log{A}}}{1+\frac{1}{\log{A}}+\frac{k}{\log^2{A}}}.$$ Since $A\geq 10^5$, we get $1\leq 1+\frac{1}{\log{A}}+\frac{k}{\log^2{A}}<1.1061,$ and $1.33486<k\left(1-\frac{1}{\log{A}}\right)<2.53816 \implies -1.53816<1-k+\frac{k}{\log{A}}<-0.33486$. The result follows.
\end{proof}

\begin{proof}[Proof of Lemma \ref{prop:piAestimates}]  
By \cite[(3.5)]{RS1962}, we have $\pi(A)>\frac{A}{\log{A}}$ if $A\geq 17$. Each inequality in \eqref{prop:4.1&2&3} then follows by checking for $3\leq A<17$. The second statement of \eqref{prop:4.4&5} follows from Lemma \ref{prop:sumiOverpi} for $A\geq 10^5$, and we use computation to verify it for smaller $A$. For the first claim in \eqref{prop:4.4&5}, we use Lemma \ref{prop:DusartPi}. If $A\geq 11813$, then $$
    \frac{A}{\pi(A)}-\log{A}+1 \leq \frac{\frac{-0.46184}{\log{A}}+\frac{1.46184}{\log^2{A}}}{1+\frac{1}{\log{A}}+\frac{1.46184}{\log^2{A}}}, $$ which is negative in this range of $A$. Again, we split the interval $[3,11813)$ into $[p_k,p_{k+1})$ for $2\leq k \leq 1415$, in which $$\frac{A}{\pi(A)}-\log{A}+1=\frac{A}{k}-\log{A}+1,$$ whose derivative is $\frac{1}{k}-\frac{1}{A}>0$. Thus $\frac{p_{k+1}}{k}-\log{p_{k+1}}+1$ is an upper bound for $\frac{A}{\pi(A)}-\log{A}+1$ on this interval. Computing this for all $k$ and combining yields the values given in Table \ref{tab:AoverPiA} and the result.
\end{proof}
\begin{rem}
    Table \ref{tab:AoverPiA} is not optimal: for instance, $(1.353,7)$ is a valid choice of $(c,J)$, but we have chosen to omit it because it provides little improvement.
\end{rem}
\begin{lemma} \label{prop:liBound}
    Suppose that $A\geq 10^5$. Then $$\frac{6A}{\log^4{A}}<\li{A}-\frac{A}{\log{A}}-\frac{A}{\log^2{A}}-\frac{2A}{\log^3{A}}<\frac{10.5A}{\log^4{A}}.$$
\end{lemma}
\begin{proof}
    This holds at $A=10^5$; then consider derivatives on $A> 10^5$.
\end{proof}
\begin{lemma}\label{prop:asympPunchConj}
    Let $c$ be a real number in the interval $(0,1)$ and let $i$ range over the positive integers. Then $\pi(cp_i)\sim ci$ and $i-\pi(cp_i)\sim (1-c)i$ as $i\rightarrow\infty$.
\end{lemma}
\begin{proof}
    By the prime number theorem, $\pi(i)\sim\frac{i}{\log{i}}$ and $p_i\sim i\log{i}$ as $i\rightarrow\infty$, so $$
        \pi(cp_i)\sim \frac{c p_i}{\log{c p_i}}\sim \frac{ci \log{i}}{\log{c}+\log{i}+\log{\log{i}}}
        \sim \frac{ci\log{i}}{\log{i}}=ci.$$ The second claim follows immediately.
\end{proof}
\begin{lemma} \label{prop:p_nLowerBound} \cite[Theorem 3]{Dusart1999} If $n\geq 2$, then $p_n\geq n(\log{n}+\log{\log{n}}-1).$
\end{lemma}
\begin{lemma} \label{prop:p_nUpperBoundBetter} \cite[Section 4]{Dusart1999} If $n\geq 39017$, then $p_n\leq n(\log{n}+\log{\log{n}}-0.9484).$
\end{lemma}
\begin{lemma} \label{prop:p_nUpperBoundBest} \cite[Section 4]{Dusart1999} If $n\geq 27076$, then $p_n\leq n(\log{n}+\log{\log{n}}-1+\frac{-1.8+\log{\log{n}}}{\log{n}})$. \end{lemma}

\begin{lemma}\label{prop:PunchWinCondition}
Fix a positive integer $i$ and suppose that there is a positive integer $k\leq i$ such that $i p_j/p_i\leq j$ for every $1\leq j \leq k$. Let $c\in [0,1]$; if $c\leq \frac{p_k}{p_i}$, then $i-\pi(cp_i)-1\leq (1-c)i$.
\end{lemma}
\begin{proof}
Observe that the result is equivalent to $\pi(cp_i)+1\geq ci$. We allow $c$ to vary over $[0,1]$. It is clear that $$\pi(cp_i)+1=\begin{cases}
    1 & 0\leq c<p_1/p_i\\
    2 & p_1/p_i\leq c < p_2/p_i \\
    \ldots \\
    i & p_{i-1}/p_i\leq c <1 \\
    i+1 & c=1.
\end{cases}$$ Since $c\leq p_k/p_i$, either it is $p_k/p_i$ or it is in one of $$[0,p_1/p_i),[p_1/p_i,p_2/p_i),\ldots,[p_{k-1}/p_i,p_k/p_i).$$ If it is in $[0,p_1/p_i)$, then we have $\pi(cp_i)+1=1$ and $ci< ip_1/p_i\leq 1$ by assumption. Similarly, if it is in $[p_{j-1}/p_i,p_j/p_i)$ for some $1<j\leq k$, then $\pi(cp_i)+1=j$ and $ci<ip_j/p_i\leq j$. Finally, if it is $p_k/p_i$, then $\pi(cp_i)+1=k+1$ but $ci=ip_k/p_i\leq k$.
\end{proof}
\begin{lemma}\label{prop:PunchConditionSatisfiedPart2}
    Fix a positive integer $i\geq 39017$, let $u=\exp(-0.0516)\approx 0.9497$ and let $k=k(i):=\floor{ui}$. Then $ip_j/p_i\leq j$ for every $1\leq j \leq k$.
\end{lemma}
\begin{proof}
    It will be more convenient to prove the statement in the following form: $i/p_i\leq j/p_j$ for every $1\leq j \leq k$. We use Lemmas \ref{prop:p_nLowerBound} and \ref{prop:p_nUpperBoundBetter}: since $i\geq 39017$, we have $$\frac{1}{\log{i}+\log{\log{i}}-0.9484}\leq \frac{i}{p_i}\leq \frac{1}{\log{i}+\log{\log{i}}-1}.$$

    Let $w=1/(\log{i}+\log{\log{i}}-1)$; then $i/p_i\leq w$. It also follows that $\frac{i\log{i}}{e}=\exp(1/w)$. We now want to prove that $j/p_j\geq w$ for every $1\leq j \leq k$.
    
    If $j\geq 39017$, it is enough to show that $\frac{1}{\log{j}+\log{\log{j}}-0.9484}\geq w$. To this end, consider that we have $j\leq k\leq \exp(-0.0516)i$, which implies that $0\leq \frac{\log{j}}{\log{i}}<1$. It follows that $\frac{j\log{j}}{i\log{i}}<\exp(-0.0516)$. This implies that $\exp(-0.9484)j\log{j}<\exp(-1)i\log{i}=\exp(1/w)$, so that $\frac{1}{\log{j}+\log{\log{j}}-0.9484}>w$ as desired.

    Suppose $j<39017$. One may check that the minimum value of $j/p_j$ in this range is greater than $0.0834618$. Further, if $i\geq 40901$, then $i/p_i\leq \frac{1}{\log{40901}+\log{\log{40901}}-1}<0.0834618$, so we are done. It remains to verify that $\min_{1\leq j \leq k(i)}\frac{j}{p_j}\geq \frac{i}{p_i}$ whenever $39017\leq i<40901$. We outline an efficient way of doing so. Let $i=39017$, then $k=37054$. It is readily verified that $\min_{1\leq j \leq k(i)}\frac{j}{p_j}=\frac{37053}{441443}>0.0839$ and that $i/p_i=39017/467473<0.0835$; the result therefore holds for this $i$. Now, consider $i+1$; we get that $k(i+1)=\floor{u(i+1)}$ is either $k(i)$ or $k(i)+1$. If it is $k(i)$, then $\min_{1\leq j \leq k(i+1)}\frac{j}{p_j}=\min_{1\leq j \leq k(i)}\frac{j}{p_j}$. If it is $k(i)+1$, then $$\min_{1\leq j \leq k(i+1)}\frac{j}{p_j}=\min\left(\min_{1\leq j \leq k(i)}\frac{j}{p_j},\frac{k(i)+1}{p_{k(i)+1}}\right);$$ this is easy to compute. Compute $i/p_i$; if it is not greater than this minimum, then we have proved the statement for $i+1$. Proceed inductively until $i=40901$.

\end{proof}

\begin{lemma}\label{prop:PunchPrimeRatioBetter}
    If $i\geq 33$ is a positive integer, $u=\exp(-0.0516)$ and $k=\floor{ui}$, then $p_k/p_i\geq 0.8759$.
\end{lemma}
\begin{proof} For brevity, we will write $y:=\frac{1}{u}$. Again, we employ Lemmas \ref{prop:p_nLowerBound} and \ref{prop:p_nUpperBoundBetter}. Firstly, the assumption $k=\floor{ui}$ implies that $ky\leq i <(k+1)y$. Fix $k$ and allow $i$ to be an integer in this interval. If we suppose that $i\geq 39017$, we will certainly have $k\geq 2$, so that \begin{align*}
        \frac{p_k}{p_i}&\geq \frac{k(\log{k}+\log{\log{k}}-1)}{i(\log{i}+\log{\log{i}}-0.9484)}> \frac{k(\log{k}+\log{\log{k}}-1)}{(k+1)y(\log{((k+1)y)}+\log{\log{((k+1)y)}}-0.9484)} \\
        &=\frac{k(\log{k}+\log{\log{k}}-1)}{y(k+1)(\log{(k+1)}+\log{(\log{(k+1)}+0.0516})-0.8968)}.
    \end{align*}
    For any positive real $x$ and $\varepsilon$, we have $\log{(x+\varepsilon)}< \frac{\varepsilon}{x}+\log{x}$. Applying this repeatedly yields $\log{(\log{(k+1)}+0.0516)}< \log{\log{k}}+\frac{1}{k\log{k}}+\frac{0.0516}{\log{k}}$ for $k>1$. Then \begin{align}\notag 
        \frac{p_k}{p_i}\geq &\frac{uk(\log{k}+\log{\log{k}}-1)}{(k+1)\left[\log{k}+\frac{1}{k}-0.8968+\log{\log{k}}+\frac{1}{k\log{k}}+\frac{0.0516}{\log{k}}\right]}\\
        \notag &=\frac{u(1-\frac{k}{k\log{k}+k\log{\log{k}}})}{1+\frac{-0.8968k+\frac{0.0516k}{\log{k}}+\log{k}+\log{\log{k}}+0.1032+\frac{1.0516}{\log{k}}+\frac{1}{k}+\frac{1}{k\log{k}}}{k\log{k}+k\log{\log{k}}}}\\
        \label{eq:primeRatioLemmaInequality}&=\frac{u\left(1-\frac{1}{\log{k}+\log{\log{k}}}\right)}{1+\frac{-0.8968+\frac{0.0516}{\log{k}}+\frac{\log{k}}{k}+\frac{\log{\log{k}}}{k}+\frac{0.1032}{k}+\frac{1.0516}{k\log{k}}+\frac{1}{k^2}+\frac{1}{k^2\log{k}}}{\log{k}+\log{\log{k}}}}.
    \end{align}

   Since $k\geq 37054$, $u\left(1-\frac{1}{\log{k}+\log{\log{k}}}\right)>0.8759$. In \eqref{eq:primeRatioLemmaInequality}, it is not hard to see that the numerator of the fraction in the denominator is decreasing and negative when $k>6$, and must therefore be in $(-0.8968,0)$. Now, $\log{k}+\log{\log{k}}$ is at least $12.87$; the denominator in \eqref{eq:primeRatioLemmaInequality} is then a positive number that is less than $1$. The entire fraction is therefore greater than $0.8759$.

    Suppose instead that $33\leq i\leq 39016$. It is easily verified by computation that $p_k/p_i\geq 0.8759$ for all such $i$.
\end{proof}

\begin{table}
    \centering
    \begin{tabular}{|c|c|c|}
    \hline
        $n$ & $i_n$ & $c_n$\\
    \hline
        1 & 5 & 48/61\\
        2 & 62 & 103/121\\
        3 & 122 & 224/263\\
    \hline
    \end{tabular}
    \caption{Values $i_n$ and $c_n$ as in Theorem \ref{prop:PunchConjecture} and Corollary \ref{prop:PunchConjectureUniform}}
    \label{tab:PunchConjectureValues}
\end{table}
\begin{theorem}\label{prop:PunchConjecture}
Let $n\in \{1,2,3\}$, and let $i_n$ and $c_n$ be as in Table \ref{tab:PunchConjectureValues}. Then for any real number $c\in[0,c_n]$ and any integer $i\geq i_n$, we have $i-\pi(cp_i)-1\leq(1-c)i$.
\end{theorem}
\begin{proof} Suppose first that $i\geq 39017$. Let $u=\exp(-0.0516)$ and $k=\floor{ui}$. By Lemma \ref{prop:PunchPrimeRatioBetter}, $p_k/p_i\geq 0.8759$, which exceeds every $c_n$. Further, by Lemma \ref{prop:PunchConditionSatisfiedPart2}, $ip_j/p_i\leq j$ for every $1\leq j\leq k$. By Lemma \ref{prop:PunchWinCondition}, these two facts imply the result for $i\geq 39017$.

When $i_n \leq i \leq 39016$, we apply Lemma \ref{prop:PunchWinCondition}. Fix an $i$, select the least $k$ such that $c_n<p_k/p_i$; then verify by computation that $i/p_i\leq j/p_j$ for every $1\leq j \leq k$. When $n=1$, this proves the result for every integer $i\in[i_n,39016]\setminus \{6,8,11,14,61\}$. When $n=2$, the exceptions are $\{121,263\}$. When $n=3$, the sole exception is $263$.

For the exceptional values of $i$, we use Lemma \ref{prop:PunchWinCondition} in a slightly different way. To illustrate the concept, we use $n=1$, $i=11$ as an example. Select the greatest $k$ such that $p_k/p_i\leq c_n=\frac{48}{61}$; here, this is $k=9$. One can verify that $ip_j/p_i\leq j$ for every $1\leq j \leq k$; so we are done if $c\leq p_k/p_i=23/31$. To complete the proof, suppose that $p_k/p_i< c\leq \frac{48}{61}$. By definition, we must have $\pi(cp_i)=k$. Thus $i-\pi(cp_i)-1=11-9-1=1$, while $(1-c)i\geq \frac{13}{61}\cdot i=2\frac{21}{61}>1$.

The proof proceeds similarly for all of the exceptional values; see Table \ref{tab:PunchConjectureVerification}.

\begin{table}
    \centering
    \def\arraystretch{1.2}
    \begin{tabular}{|c|c|c|c|c|c|}
    \hline
         $n$ & $i$ & $k$ & $i-\pi(c_np_i)-1$ & $(1-c_n)i$\\ \hline
         1 & 6 & 4 & 1 & $1\frac{17}{61}$\\
         1 & 8 & 6 & 1 & $1\frac{43}{61}$\\  
         1 & 11 & 9 & 1 & $2\frac{21}{61}$\\  
         1 & 14 & 11 & 2 & $2\frac{60}{61}$\\  
         1 & 61 & 47 & 13 & 13\\
         2 & 121 & 102 & 18 & 18 \\
         2 & 263 & 223 & 39 & $39\frac{15}{121}$ \\
         3 & 263 & 223 & 39 & 39 \\ \hline
    \end{tabular}
    \caption{Verifying the remaining cases}
    \label{tab:PunchConjectureVerification}
\end{table}
\end{proof}

\begin{rem}
    In fact, one may select any real number $\sigma\in[0,1)$ and find some $M(\sigma)\geq 1$ such that $i-\pi(cp_i)-1\leq (1-c)i$ whenever $i\geq M(\sigma)$ and $c\in[0,\sigma]$. When $\sigma=0$, we can trivially take $M(0)=1$; otherwise, this claim can be proved by using the prime number theorem in the form $\pi(x)=\li(x)+O\left(xe^{-a\sqrt{\log{x}}}\right).$ This approach could also be used to give an explicit value for $M(\sigma)$, but the author believes that the argument of Theorem \ref{prop:PunchConjecture} is far more suitable. We note here that without using something stronger than Lemma \ref{prop:p_nUpperBoundBetter}, that argument cannot prove the existence of $M(\sigma)$ for $\sigma\geq \exp(-0.0516)$. For these $\sigma$, we would need a statement of the form $p_n\leq n(\log{n}+\log{\log{n}}-1+o(1)),$ such as Lemma \ref{prop:p_nUpperBoundBest}. 
\end{rem}
\begin{cor} \label{prop:PunchConjectureUniform}
        Select $n\in \{1,2,3\}$ and let $i_n$ and $c_n$ be as in Table \ref{tab:PunchConjectureValues}. Then for any real number $c\in[0,1]$ and any integer $i\geq i_n$, we have $i-\pi(cp_i)-1\leq i\cdot \max\left(1-c,1-c_n\right).$
\end{cor}
\begin{proof}
    Let $i\geq i_n$. If $c\in[0,c_n]$, we are done. If $c\in (c_n,1]$, then $i-\pi(c p_i)-1\leq i-\pi(c_n p_i)-1\leq (1-c_n)i.$
\end{proof}
\begin{lemma} \label{prop:EiEstimate}
Let $x\geq 17.26$. Then $\Ei(x)=\frac{e^x}{x}+\frac{e^x}{x^2}+\frac{ke^x}{x^3},$ for some $k\in (2,2.47).$
\end{lemma}
\begin{proof}
    Let $F(x)=\left(\Ei(x)-\frac{e^x}{x}-\frac{e^x}{x^2}\right)\frac{x^3}{e^x}=:uv$, say, and note that $F(17.26)<2.47$. Observing that $\frac{du}{dx}=\frac{2}{v}$ and $\frac{dv}{dx}=v\left(\frac{3}{x}-1\right)$, we see that $F$ satisfies the differential equation $F'(x)=2+F(x)\left(\frac{3}{x}-1\right)=-\left(1-\frac{3}{x}\right)\left(F(x)-\frac{2x}{x-3}\right)$. For the upper bound on $k$, it is enough to prove that $F'(x)<0$ when $x\geq 17.26$, for which we just need $G(x):=F(x)-\frac{2x}{x-3}>0$. Substituting, we get $G'(x)+\left(1-\frac{3}{x}\right)G=\frac{6}{(x-3)^2}$, and via an integrating factor we reach $\frac{d}{dx}\left[\frac{e^x}{x^3}G(x)\right]=\frac{6e^x}{x^3(x-3)^2}$. The right-hand side is positive whenever $x>3$, and by calculation we establish that $G(17.26)>0$. It follows that $G(x)>0$ when $x\geq 17.26$. For the lower bound, observe that $H(x):=\Ei(x)-\frac{e^x}{x}-\frac{e^x}{x^2}-\frac{2e^x}{x^3}$ has positive derivative for $x>0$ and that $H(17.26)>0$.
\end{proof}

\begin{lemma} \label{prop:gEstimate}
    Select $r\in\{2,3,\ldots,15\}$ and let $I_r$ be as in Table \ref{tab:gValues}. Let $$g=g(A,r):=\frac{\sum_{\substack{a\leq A\\a \text{ prime}}}a^{1-1/r}}{\pi(A)A^{1-1/r}/(2-1/r)}.$$ If $A\geq 10^5$, then $g=1-\frac{\frac{r-1}{2r-1}}{\log{A}}+\frac{\gamma(r)}{\log^2{A}}$ for some $\gamma(r)\in I_r$.
\end{lemma}
\begin{proof} By partial summation, $$g=\left(2-\frac{1}{r}\right)-\left(1-\frac{1}{r}\right)\left(2-\frac{1}{r}\right)\cdot \frac{\int_{2}^{A}\pi(t)t^{-1/r}\,dt}{\pi(A)A^{1-1/r}}.$$ We split the integral into $\int_{2}^{11813}$ and $\int_{11813}^{A}$, and use Lemma \ref{prop:DusartPi} to deal with $\pi(t)$. This leads us to the following integrals, which we treat using Lemma \ref{prop:EiEstimate} (permissible because $\left(2-\frac{1}{r}\right)\log{A}\geq \frac{3}{2}\log{10^5}>17.26$). We have $$\int_{11813}^{A}\frac{t^{1-1/r}}{\log{t}}\,dt=\frac{A^{2-1/r}}{(2-1/r)\log{A}}+\frac{A^{2-1/r}}{(2-1/r)^2\log^2{A}}+c_1\frac{A^{2-1/r}}{(2-1/r)^3\log^3{A}}+K_1,$$ where $K_1=-\Ei\left[(2-1/r)\log{11813}\right]$ and $c_1\in (2,2.47)$,
$$\int_{11813}^{A}\frac{t^{1-1/r}}{\log^2{t}}=\frac{A^{2-1/r}}{(2-1/r)\log^2{A}}+c_2 \cdot \frac{A^{2-1/r}}{(2-1/r)^2\log^3{A}}+K_2,$$ where $K_2=-(2-1/r)\Ei\left[(2-1/r)\log{11813}\right]+\frac{11813^{2-1/r}}{\log{11813}}$ and $c_2\in (2,2.47)$, and 
$$\int_{11813}^{A}\frac{t^{1-1/r}}{\log^3{t}}\,dt=\frac{c_3}{2}\cdot \frac{A^{2-1/r}}{(2-1/r)\log^3{A}}+K_3,$$ where $K_3=-\frac{1}{2}\left[(2-1/r)^2\Ei\left[(2-1/r)\log{11813}\right]-\frac{11813^{2-1/r}\left(1+(2-1/r)\log{11813}\right)}{\log^2{11813}}\right]$ and $c_3\in (2,2.47)$. Write $$K_4=\int_{2}^{11813}\pi(t)t^{-1/r}\,dt=\frac{r}{r-1}\left(1416p_{1416}^{1-1/r}-\sum_{i=1}^{1416}p_i^{1-1/r}\right);$$ the original integral becomes \begin{align*}&K_4+\int_{11813}^{A}\frac{t^{1-1/r}}{\log{t}}+\frac{t^{1-1/r}}{\log^2{t}}+c_4\cdot \frac{t^{1-1/r}}{\log^3{t}}\,dt \\
&=A^{2-1/r}\left[\frac{\frac{1}{2-1/r}}{\log{A}}+\frac{\frac{1}{2-1/r}+\frac{1}{(2-1/r)^2}}{\log^2{A}}+c_5\frac{\frac{1}{2-1/r}}{\log^3{A}}\right]+K_5,
\end{align*} where $c_4\in[1.46184,2.53816]$ and we have written $c_5=\left(\frac{c_1}{(2-1/r)^2}+\frac{c_2}{2-1/r}+\frac{c_3 c_4}{2}\right)$ and $K_5=(K_1+K_2+K_4+c_4 K_3)$. We get $$g=\left(2-\frac{1}{r}\right)-\left(1-\frac{1}{r}\right) \frac{A}{\pi(A)}\left(\frac{1}{\log{A}}+\frac{3-1/r}{(2-1/r)\log^2{A}}+\frac{c_5}{\log^3{A}}+\frac{K_5(2-1/r)}{A^{2-1/r}}\right).$$ Via Lemma \ref{prop:AoverPiA}, $\frac{A}{\pi(A)}=\log{A}-1+\frac{c_6}{\log{A}}$ where $c_6\in (-1.53816,-0.30273)$. Then \begin{align*}g&=1-\frac{\frac{r-1}{2r-1}}{\log{A}}-\frac{1-1/r}{\log^2{A}}\biggl[c_5+c_6-\frac{3-1/r}{2-1/r}+\frac{c_6\cdot \frac{3-1/r}{2-1/r}-c_5}{\log{A}}+\frac{c_5 c_6}{\log^2{A}}\\&+\frac{K_5(2-1/r)\log^3{A}}{A^{2-1/r}}-\frac{K_5(2-1/r)\log^2{A}}{A^{2-1/r}}+\frac{c_6K_5(2-1/r)\log{A}}{A^{2-1/r}}\biggr].
\end{align*} The result follows by estimating the size of the bracketed term, based on the ranges of possible values for $c_5,c_6$ and $K_5$, and the assumption $A\geq 10^5$.
\end{proof}

    \begin{table}
        \centering
        \begin{tabular}{|cc|cc|}
        \hline
           $r$  & $I_r$ & $r$ & $I_r$\\
        \hline
           $2$ & $(-1.96259,0.21192)$ & $9$ & $(-2.94510,0.65145)$\\
           $3$ & $(-2.41178,0.38444)$ & $10$ & $(-2.97038,0.66560)$\\
           $4$ & $(-2.61876,0.48055)$ & $11$ & $(-2.99098,0.67723)$\\
           $5$ & $(-2.73863,0.54068)$ & $12$ & $(-3.00810,0.68696)$\\
           $6$ & $(-2.81698,0.58168)$ & $13$ & $(-3.02254,0.69520)$\\
           $7$ & $(-2.87225,0.61137)$ & $14$ & $(-3.03490,0.70229)$\\
           $8$ & $(-2.91334,0.63384)$ & $15$ & $(-3.04559,0.70845)$\\
        \hline
        \end{tabular}
        \caption{Intervals $I_r$ as in Lemma \ref{prop:gEstimate}}
        \label{tab:gValues}
    \end{table}

\begin{lemma} \label{prop:BinomialG}
Select $r\geq 2$. For any real $h$ with $\lvert h\rvert\leq \varepsilon$ for some $\varepsilon<1$, we have \begin{align}\label{eq:BinomialG}\left(\frac{2r-2}{2r-2-h}\right)^{2r}-1=\frac{r}{r-1}h+\gamma(r,\varepsilon)h^2,\end{align} where $\gamma(r,\varepsilon)\in \left[\frac{e^{-\frac{r}{r-1}\varepsilon}+\frac{r}{r-1}\varepsilon-1}{\varepsilon^2},\frac{(1-\varepsilon)^{-\frac{r}{r-1}}-\frac{r}{r-1}\varepsilon-1}{\varepsilon^2}\right]$.
\end{lemma}
\begin{proof} One can show that $1-h\leq \left(1-\frac{h}{k}\right)^k\leq e^{-h}$ for any $\lvert h\rvert <1$ and $k>1$. Now, the left-hand side is $\left(1-\frac{h}{2r-2}\right)^{(2r-2)\cdot -\frac{r}{r-1}}-1,$ which must therefore be in the interval $\left[e^{\frac{r}{r-1}h}-1,(1-h)^{-\frac{r}{r-1}}-1\right].$ Thus $$\frac{e^{\frac{r}{r-1}h}-\frac{r}{r-1}h-1}{h^2}\leq\frac{\left(\frac{2r-2}{2r-2-h}\right)^{2r}-\frac{r}{r-1}h-1}{h^2}\leq \frac{(1-h)^{-\frac{r}{r-1}}-\frac{r}{r-1}h-1}{h^2}.$$ Consider the left-hand side: letting $\omega=\frac{r}{r-1}h$, we have $\left(\frac{r}{r-1}\right)^2\frac{e^\omega-\omega-1}{\omega^2}$. Since $\frac{e^\omega-\omega-1}{\omega^2}$ is increasing for all real $\omega$, this is at least $\frac{e^{-\frac{r}{r-1}\varepsilon}+\frac{r}{r-1}\varepsilon-1}{\varepsilon^2}.$ Similarly, for any $c>1$, the function $\frac{(1-x)^{-c}-cx-1}{x^2}$ is increasing on $x<1$, so the right-hand side is no more than $\frac{(1-\varepsilon)^{-\frac{r}{r-1}}-\frac{r}{r-1}\varepsilon-1}{\varepsilon^2}.$
\end{proof}

    \begin{table}
        \centering
        \begin{tabular}{|cc|cc|}
        \hline
           $r$  & $I_r$ & $r$ & $I_r$\\
        \hline
           $2$ & $(-3.74838,1.23227)$ & $9$ & $(-3.24448,1.40792)$\\
           $3$ & $(-3.49571,1.31890)$ & $10$ & $(-3.23389,1.41173)$\\
           $4$ & $(-3.39187,1.35538)$ & $11$ & $(-3.22533,1.41481)$\\
           $5$ & $(-3.33543,1.37540)$ & $12$ & $(-3.21828,1.41735)$\\
           $6$ & $(-3.29999,1.38803)$ & $13$ & $(-3.21237,1.41948)$\\
           $7$ & $(-3.27568,1.39673)$ & $14$ & $(-3.20735,1.42129)$\\
           $8$ & $(-3.25796,1.40308)$ & $15$ & $(-3.20303,1.42285)$\\
        \hline
        \end{tabular}
        \caption{Intervals $I_r$ as in Lemma \ref{prop:eq47Part2}}
        \label{tab:eq47values}
    \end{table}
\begin{lemma}\label{prop:eq47Part2}
Select $r\in \{2,\ldots, 15\}$. Let $I_r$ be as in Table \ref{tab:eq47values} and $g$ as in Lemma \ref{prop:gEstimate}. If $A\geq 10^5$, then $\left(\frac{2r-2}{2r-1-g}\right)^{2r}=1-\frac{\frac{r}{2r-1}}{\log{A}}+\frac{\gamma(r)}{\log^2{A}}$ for some $\gamma(r)\in I_r$.
\end{lemma}
\begin{proof}
    By Lemma \ref{prop:gEstimate}, we get (by verifying that $A$ is sufficiently large) that $$\lvert g-1 \rvert=\frac{\frac{r-1}{2r-1}}{\log{A}}-\frac{\gamma_1(r)}{\log^2{A}}\leq \frac{\frac{r-1}{2r-1}}{\log{10^5}}-\frac{\gamma_1(r)}{\log^2{10^5}},$$ where $\gamma_1(r)$ is as in Table \ref{tab:gValues}. Take the maximum possible value of the right-hand side (using the lower bound on $\gamma_1(r)$) as our $\varepsilon$ in Lemma \ref{prop:BinomialG} with $h=g-1$; this leads to $$\left(\frac{2r-2}{2r-1-g}\right)^{2r}=1+\frac{r}{r-1}(g-1)+\gamma_2(r)(g-1)^2,$$ with $\gamma_2(r)$ in a certain interval, which we calculate but will not write down. Thus \begin{align*}
        \left(\frac{2r-2}{2r-1-g}\right)^{2r}&=1-\frac{\frac{r}{2r-1}}{\log{A}}+\\&\frac{1}{\log^2{A}}\left[\gamma_1(r)\frac{r}{r-1}+\gamma_2(r)\left(\frac{r-1}{2r-1}\right)^2-\frac{2\gamma_1(r)\gamma_2(r)\frac{r-1}{2r-1}}{\log{A}}+\frac{\gamma_2(r)\gamma_1(r)^2}{\log^2{A}}\right].
    \end{align*} The result follows from $A\geq 10^5$ and the bounds on $\gamma_1(r)$ and $\gamma_2(r)$.
\end{proof}
\begin{lemma}\label{prop:sumiOverpi}
    Suppose that $A\geq 10^5$. Then $-1+\frac{A}{\pi(A)^2}\sum_{\substack{a\leq A\\a\text{ prime}}}\frac{\pi(a)}{a}=\frac{\ell}{\log{A}},$ where $\ell\in (0.933,1.874)$.
\end{lemma}
\begin{proof}
    By Lemma \ref{prop:DusartPi}, we have $$\sum_{\substack{a\leq A\\ a\text{ prime}}}\frac{\pi(a)}{a}<1051.43+\sum_{\substack{10^5<a\leq A\\ a\text{ prime}}}\frac{\pi(a)}{a}\leq 1051.43+\sum_{\substack{10^5<a\leq A\\ a\text{ prime}}}\left(\frac{1}{\log{a}}+\frac{1.2762}{\log^2{a}}\right).$$ Via partial summation, the sum becomes $$
        \pi(A)\left(\frac{1}{\log{A}}+\frac{1.2762}{\log^2{A}}\right)-\pi(10^5)\left(\frac{1}{\log{10^5}}+\frac{1.2762}{\log^2{10^5}}\right)+\int_{10^5}^{A}\pi(t)\left(\frac{1}{t\log^2{t}}+\frac{2.5524}{t\log^3{t}}\right)\,dt.$$ Isolate the integral, then apply Lemma \ref{prop:DusartPi} and Lemma \ref{prop:liBound}: \begin{align*}
            &\int_{10^5}^{A}\left(\frac{1}{\log{t}}+\frac{1.2762}{\log^2{t}}\right)\left(\frac{1}{\log^2{t}}+\frac{2.5524}{\log^3{t}}\right)\,dt\leq \int_{10^5}^{A}\left(\frac{1}{\log^3{t}}+\frac{3.8286}{\log^4{t}}+\frac{3.2574}{\log^5{t}}\right)\,dt \\
            &< 1.273825\li{A}-\frac{1.273825A}{\log{A}}-\frac{1.273825A}{\log^2{A}}-\frac{1.54765A}{\log^3{A}}-\frac{0.81435A}{\log^4{A}}-135.3\\&<\frac{A}{\log^3{A}}+\frac{12.6A}{\log^4{A}}-135.3.
        \end{align*} Combine and apply Lemma \ref{prop:AoverPiA}; we get \begin{align*}
            -1+\frac{A}{\pi(A)^2}\sum_{\substack{a\leq A\\ a\text{ prime}}}\frac{\pi(a)}{a}&<\frac{A}{\pi(A)}\left(\frac{1}{\log{A}}+\frac{1.2762}{\log^2{A}}\right)+\frac{A^2}{\pi(A)^2\log^3{A}}+\frac{12.6A^2}{\pi(A)^2\log^4{A}}-1\\
            &<\frac{1.2762}{\log{A}}+\frac{9.0211}{\log^2{A}}-\frac{25.1918}{\log^3{A}}+\frac{5.5767}{\log^4{A}}+\frac{7.7205}{\log^5{A}}+\frac{1.1548}{\log^6{A}}.
        \end{align*} From this, it is not difficult to show that $$\left(-1+\frac{A}{\pi(A)^2}\sum_{\substack{a\leq A\\a\text{ prime}}}\frac{\pi(a)}{a}\right)\cdot \log{A}<1.874.$$ The lower bound can be established in the same way.
\end{proof}
\begin{lemma}\label{prop:Weil} Let $\chi$ be the quadratic character modulo an odd prime $p$ and choose an integer $r\geq 1$. Select $c_1,\ldots,c_{2r}$ that are all distinct modulo $p$. Then $$\left\lvert\sum_{x=1}^{p}\chi\left((x+c_1)(x+c_2)\ldots (x+c_{2r})\right)\right\rvert\leq 2(r-1)\sqrt{p}+1.$$
\end{lemma}
\begin{proof}
    When $r\geq 2$, a bound of $(2r-1)\sqrt{p}$ follows from \cite[Lemma 9.25]{MV2007} by taking $d=2$, $k=2r$ and $e_j=1$ for each $j$. It can be improved to $2(r-1)\sqrt{p}+1$ because the character has order $2$. Note that the $+1$ term does not appear in \cite{Booker2006}; the author thanks Andy Booker for pointing out that this is an error.

    When $r=1$, the sum is always $-1$. To see this, note that without loss of generality we may choose $c_2=0$; we get \begin{align*}
        \sum_{x\not \equiv 0,-c_1}\chi(x(x+c_1))&=\sum_{x\not \equiv 0,-c_1}\chi\left(x^2\left(1+\frac{c_1}{x}\right)\right)=\sum_{x\not \equiv 0,-c_1}\chi\left(1+\frac{c_1}{x}\right)\\
        &=\sum_{y\not\equiv -1,0}\chi(1+y)=\sum_{y\not\equiv 0,1}\chi(y)=-1.
    \end{align*}
\end{proof}

\section*{Acknowledgements}
The author is very grateful to Tim Trudgian for his supervision, encouragement, support, feedback and patience. The author thanks Daniel Johnston for discussions about estimates for $\pi(x)$, and Andy Booker for discussions about \cite{Booker2006}. The author used GPT-6 Luna to find a way to prove Lemma \ref{prop:EiEstimate}, but the proof text itself was written by the author. Finally, the author acknowledges the receipt of an Australian Government Research Training Program Scholarship.

\bibliographystyle{plain}
\bibliography{refs}
\end{document}